\documentclass[11pt]{article}

\RequirePackage{amsthm}
\RequirePackage{amssymb}
\usepackage[all,ps]{xy}

\usepackage{latexsym}
\usepackage[frenchb]{babel}

\usepackage[latin1]{inputenc}

\usepackage{a4wide}

\author{Andrea Surroca Ortiz\\
\'Ecole Polytechnique F\'ed\'erale de Lausanne}
\title{\Large Sur l'effectivit\'e du th\'eor\`eme de Siegel et la conjecture $abc$}
\date{26 octobre 2005}

\newtheorem{thm}{Th\'eor\`eme}[section]
\newtheorem{conj}[thm]{Conjecture}
\newtheorem{hyp}[thm]{Hypoth\`ese}
\newtheorem{prop}[thm]{Proposition}
\newtheorem{lemme}[thm]{Lemme}
\newtheorem{cor}[thm]{Corollaire}
\newtheorem{remarque}[thm]{Remarque}
\newtheorem{definition}[thm]{D\'efinition}
\newtheorem{propriete}[thm]{Propri\'et\'e}

\newcommand\enteros{\mathbf{Z}}
\newcommand\rat{\mathbf{Q}}
\newcommand\Qbarre{\overline{\mathbf{Q}} }
\newcommand\espproj{{\mathrm{I} \hspace{-0,1 cm} \mathrm{P}}}
\newcommand\droitemoins{\espproj^{1}\setminus\{0,1,\infty\}}
\newcommand\Ebarre{\overline{E}}
\newcommand\card{\mathrm{card}}
\newcommand\rad{\mathrm{rad}}
\newcommand\supp{\mathrm{supp}}
\newcommand\ord{\mathrm{ord}}
\newcommand\Spec{\mathrm{Spec}}
\newcommand\Div{\mathrm{Div}}

\newcommand\esp{\hspace{0,2cm}}
\newcommand\espa{\hspace{0,17cm}}
\newcommand\Belyi{Bely\u\i}
\newcommand\PK{\mathcal{P}_K}
\newcommand\vp{\mathfrak{p}}
\newcommand\vq{\mathfrak{q}}

\newcommand\g{``}

\newcommand\vers{\longrightarrow }

\begin{document}

\maketitle

\begin{quote}
\textbf{R\'esum\'e.}
{\small Nous montrons qu'un raffinement du th\'eor\`eme de Siegel sur les points entiers de courbes alg\'ebriques impliquerait la conjecture $abc$ de Masser-Oesterl\'e. 
Nous formulons une hypoth\`ese ``Siegel uniforme"qui est une majoration de la hauteur des points $S$-entiers de la courbe, en termes du corps de rationalit\'e et de l'ensemble de places $S$. La validit\'e de l'hypoth\`ese pour une quelconque courbe alg\'ebrique de caract\'eristique d'Euler-Poincar\'e strictement n\'egative, impliquerait une version de la conjecture $abc$.
Ceci \'etend aux points $S$-entiers des r\'esultats pr\'ec\'edents de L. Moret-Bailly \cite{moret-bailly}, et est en quelque sorte, un \'enonc\'e r\'eciproque de ce que nous avons montr\'e dans \cite{parma}, en suivant les id\'ees propos\'ees par N. Elkies \cite{elkies}. Le principal outil g\'eom\'etrique employ\'e est un th\'eor\`eme de G.V. \Belyi. Nous montrons aussi quelques versions inconditionnelles de ces \'enonc\'es : un r\'esultat allant dans le sens de la conjecture $abc$, valable sur tout corps de nombres, ainsi que des bornes pour la hauteur des solutions en $S$-entiers de certaines \'equations diophantiennes classiques.

 }

\end{quote}

\begin{quote}
\textbf{Mots clefs : }
{\small  Th\'eor\`eme de Siegel, conjecture abc, fonction de \Belyi, bornes uniformes de points S-entiers.}
\end{quote}


%

\section{Introduction.}\label{intro.siegel}

\'Etant donn\'ee une courbe alg\'ebrique affine $U$ d\'efinie sur un corps de nombres $K$, de genre $g$ et ayant $t$ points ``\`a l'infini'', C.L. Siegel d\'emontra \cite{siegel} que {\it si la caract\'eristique d'Euler-Poincar\'e $\chi(U) = 2-2g-t$ est strictement n\'egative, alors la courbe $U$ n'a qu'un nombre fini de points entiers.} Gr\^ace \`a K. Mahler, le th\'eor\`eme a \'et\'e \'etendu aux points $S$-entiers.

La d\'emonstration du th\'eor\`eme de Siegel, tout comme celle du th\'eor\`eme de Faltings (conjecture de Mordell), ne fournit pas un moyen de trouver les points et, dans le cas g\'en\'eral, ceci demeure un probl\`eme ouvert. 
Dans le cas des points rationnels, des liens ont \'et\'e \'etablis entre la question de l'effectivit\'e et la c\'el\'ebre conjecture $abc$ de D. Masser et J. Oesterl\'e (\cite{masser} et \cite{oesterle}). 
Notons $h_K$ la hauteur (logarithmique) de Weil et $\rad_K$ le radical (ou support) sur l'espace $\espproj^2(K)$. (Des notations plus d\'etaill\'ees se trouvent dans la section \ref{notations}.) 

\medskip

\begin{conj}\label{abc} {$\mathbf{(abc)}$}
 \'Etant donn\'e un corps de nombres $K$ et un nombre r\'eel $\varepsilon >0$, il existe une constante $c_{\varepsilon, K} > 0$ telle que, pour tout triplet $(a,b,c)$ de nombres alg\'ebriques dans $K$, non nuls, tels que $a+b=c$, on ait
$$h_K(a:b:c) < (1+\varepsilon)\, \rad_K(a:b:c) + c_{\varepsilon,K}.$$
\end{conj}

L. Moret-Bailly \cite{moret-bailly} montra  qu'une hypoth\`ese du type ``Mordell effectif" sur la courbe d'\'equation $y^2 +y = x^5$ impliquerait une version de la conjecture $abc$. 
Par ailleurs,  N. Elkies \cite{elkies} montra comment il serait possible de majorer la hauteur des points rationnels d'une courbe de genre sup\'erieur ou \'egal \`a 2, en supposant vraie la conjecture \ref{abc}. Les outils employ\'es par N. Elkies  sont essentiellement la th\'eorie des hauteurs et un th\'eor\`eme de G. V. \Belyi \esp \cite{belyi}.

\begin{thm}{\bf (\Belyi)}\label{th.belyi}
Une courbe alg\'ebrique projective $C$ est d\'efinie sur $\Qbarre$ si et seulement s'il existe un morphisme fini et surjectif $f : C \longrightarrow \espproj^1$ non ramifi\'e en dehors de $\{0, 1, \infty\}$. De plus, \'etant donn\'e un ensemble $\Sigma$ de points alg\'ebriques de $C$, on peut choisir $f$ tel que $f(\Sigma) \subset \{0, 1, \infty\}$. 

\end{thm}

{\small La construction d'une telle fonction de \Belyi \espa est totalement explicite et montre que, si $C$ est d\'efinie sur un corps de nombres $K$, alors la fonction $f$ peut \^etre choisie d\'efinie sur $K$.}

\medskip

En suivant les id\'ees de \cite{elkies}, nous avons montr\'e,  dans un travail ant\'erieur \cite{parma}, comment la conjecture \ref{abc} impliquerait une version effective du th\'eor\`eme de Siegel. Il s'agit d'une majoration de la hauteur des points $S$-entiers des courbes alg\'ebriques de caract\'eristique d'Euler-Poincar\'e strictement n\'egative, en termes du corps de rationalit\'e et de l'ensemble de places $S$. 
 R\'eciproquement, nous formulons dans ce travail une hypoth\`ese du type \g Siegel uniforme", o\`u la borne de la hauteur des points $S$-entiers serait uniforme en l'ensemble de places $S$. L'analogue de l'hypoth\`ese \ref{se(U,K)} ci-dessous, dans la th\'eorie des corps de fonctions, est un th\'eor\`eme (cf. la section \ref{siegel-abc} pour des r\'ef\'erences).

Fixons une courbe $U$ d\'efinie sur un corps de nombres $K$ et v\'erifiant les hypoth\`eses du th\'eor\`eme de Siegel, et une fonction hauteur $h_{U}$ d\'efinie sur $U$ et associ\'ee \`a un diviseur de degr\'e 1.

\begin{hyp}{{\bf \g Siegel Uniforme" $\mathbf{(U,K)}$}}\label{se(U,K)} 

Pour tout entier naturel $\delta \geq1$, il existe des r\'eels $k_1(U,h_{U},\delta) >1, \esp k_2(U,h_{U},\delta) >1 $ et $k_3(U,h_{U},\delta) >0$ tels que, pour toute extension finie $L/K$ de degr\'e $[L: K] \leq \delta$, pour tout sous-ensemble fini $S$ de places de $L$ et pour tout point $S$-entier $x$ de $U$, on a
$$h_{U,L}(x) \leq k_1 \sum_{\vp \in S}  \log N_{L/\rat}(\vp) + k_2 \log D_L + k_3,$$
o\`u $D_L$ est la valeur absolue du discriminant du corps $L$, et $h_{U,L} = [L:\rat] h_U$.  
\end{hyp}


Dans la section \ref{siegel-abc} nous montrons comment l'hypoth\`ese \ref{se(U,K)} impliquerait une forme de la conjecture $abc$, ce qui \'etent aux points $S$-entiers et \`a toute courbe de caract\'eristique d'Euler-Poincar\'e strictement n\'egative le r\'esultat de L. Moret-Bailly.


\begin{thm}\label{SE-abc}

Soient $K$ un corps de nombres, $U$ une courbe alg\'ebrique affine d\'efinie sur $K$ telle que $\chi (U) < 0$  et $h_U$ une hauteur d\'efinie sur $U$ associ\'ee \`a un diviseur de degr\'e 1. 

Supposons v\'erifi\'ee l'hypoth\`ese \ref{se(U,K)} pour la courbe $U$.

Pour tout triplet $(a,b,c)$ d'\'el\'ements non nuls du corps $K$ v\'erifiant $a+b=c$, on a 
$$h_{K}(a:b:c) \leq  \eta_1 \rad_{K}(a:b:c) + \eta_2 \log D_K + \eta_3,$$
o\`u le nombre $\eta_2$ ne d\'epend que du degr\'e $d$ de la fonction de \Belyi \esp associ\'ee \`a la courbe $U$, et les nombres $\eta_1$ et $\eta_3$ d\'ependent en plus de $d$ et de $[K:\rat]$, de la courbe $U$.

\end{thm}


 Dans le cas particulier, mais crucial, de la droite projective priv\'ee de trois points, nous montrons  les {\bf th\'eor\`emes \ref {abc-S-unites}} et {\bf \ref{S-unites-siegelP1}}, qui  donnent l'\'equivalence entre 
 
 \indent i)  une version de la conjecture $abc$, 
 
 ii) une majoration de la hauteur des points $S$-entiers de $U= \droitemoins$, valable pour tout ensemble $S$, et 
 
 iii) une majoration de la hauteur des $S$-unit\'es solutions de l'\'equation $u+v=1$, valable pour tout ensemble $S$.

\smallskip

Nous obtenons aussi des versions faibles, mais inconditionnelles de ces r\'esultats. Premi\`e-\
rement, un \'enonc\'e allant dans le sens de la conjecture $abc$, valable sur tout corps de nombres. Un tel r\'esultat n'\'etait connu que pour $K=\rat$. Cf. \cite{stewart-tijdeman} et \cite{stewart-yu}.

\begin{thm}\label{abc.exp}

Pour tout corps de nombres $K$, il existe des nombres r\'eels $\gamma_{1}, \gamma_{2}$ et $\gamma_{3}$  effectivement calculables, tels que, pour tout triplet $(a,b,c)$ de nombres alg\'ebriques non nuls appartenant \`a $K$ et tels que $a+b=c$, on ait
$$h_K(a:b:c) \leq \exp \{\gamma_{1} \, \rad_K(a:b:c)  + \gamma_{2}\, \log D_{K} + \gamma_{3}\}.$$

\end{thm}

 On peut choisir $\gamma_1$ et $\gamma_3$ de la forme $\alpha [K:\rat] $, avec $\alpha$  constante absolue effectivement calculable et $\gamma_{2}= 2\,[K:\rat] -1$. On peut aussi choisir $\gamma_{2}$ ind\'ependante de $[K:\rat]$ (cf. la remarque \ref{bilu-abc.exp}).

Nous en d\'eduisons ensuite  une majoration de la hauteur des points $S$-entiers de $\droitemoins $ ({\bf corollaire \ref{cor.P1}}) et des majorations de la hauteur des solutions $S$-enti\`eres de plusieurs \'equations diophantiennes classiques (cf. la section \ref{eq.dio.classiques}).


\bigskip

Ce travail est organis\'e de la fa\c con suivante. Dans la section \ref{notations} nous fixons les notations utilis\'ees concernant les corps de nombres, d\'efinissons quelques notions de hauteur, ainsi que le radical. Nous y d\'emontrons quelques lemmes pr\'eliminaires, dont une version affine du th\'eor\`eme de Chevalley-Weil et un lemme concernant le discriminant d'une extension de corps de nombres, qui nous serons utiles \`a plusieurs reprises. 

Dans la section \ref{siegel-abc}, nous expliquons les motivations \`a l'hypoth\`ese \ref{se(U,K)} et donnons la d\'emons-\\
tration du th\'eor\`eme \ref{SE-abc}. L'\'etude du cas particulier o\`u la courbe est la droite projective priv\'ee de 0, 1 et l'infini se trouve dans la section \ref{droite.proj}. Dans la section
\ref{meth.baker}, nous  rendons effectifs, en termes de hauteur, les \'enonc\'es permettant de d\'eduire la finitude des points entiers d'une courbe \`a partir de celle d'une autre, si certains morphismes lient les deux courbes.

La derni\`ere section est consacr\'ee \`a des applications des r\'esultats pr\'ec\'edents. Dans le premier paragraphe de la section \ref{applications}, se trouve la d\'emonstration du th\'eor\`eme \ref{abc.exp}. Au paragraphe \ref{genre.nul} nous \'enon\c cons le r\'esultat sur les courbes de genre nul.  On \'enonce ensuite, au paragraphe \ref{eq.dio.classiques}, les r\'esultats concernant d'autres courbes.


\medskip


\section{Notations et pr\'eliminaires.}\label{notations}


Dans tout le texte, $K$ d\'esigne un corps de nombres, $D_{K}$ la valeur absolue de son discriminant et $M_{K}$ l'ensemble de classes d'\'equivalence de ces valeurs absolues normalis\'ees de fa\c con \`a ce que la formule du produit pour un \'el\'ement non nul $x$ de $K$ s'\'ecrive $\prod_{v \in M_{K}} |x|_v^{d_{v}}=1$, o\`u $d_{v}$ est le degr\'e local en la place $v$.  

Notons $\mathcal{P}_{K}$ l'ensemble des id\'eaux premiers de $K$ et $\mathcal{P}$ l'ensemble des nombres premiers. Si $\vp \in \mathcal{P}_{K}$ \textit{divise} le nombre premier $p$ (i.e. $\vp | p$), on note $e_{\vp}= \ord_{\vp}(p)$ l'indice de ramification de $\vp$ au-dessus de $p$ et $f_{\vp}$ le degr\'e du corps r\'esiduel. Alors 
\begin{equation}
d_{\vp} = [K_{\vp}:\rat_p]= e_{\vp} f_{\vp}, \esp  \sum_{\vp |p} e_{\vp} f_{\vp} = [K:\rat] \esp \textrm{et} \esp N_{K/\rat}(\vp) = p^{f_{\vp}} .  \label{degre-norme}
\end{equation}
 On pose $|x|_{\vp} = p^{-\ord_{\vp}(x)/e_{\vp}}$ et, pour toute place $v$ de $M_K$, on pose $v(x) = - \log |x|_v$, avec la convention $v(0)=\infty$.

Si $S$ est un ensemble fini de places du corps $K$ contenant les places archim\'ediennes, $O_{K}$ et $O_{K,S}$ d\'esignent respectivement l'anneau des entiers et l'anneau des $S$-entiers de $K$. 

Si $v$ est une place ultram\'etrique de $K$ correspondant \`a l'id\'eal premier $\vp$, $N(v)=N_{K/\rat}(\vp)$ d\'esigne la norme de l'id\'eal $\vp$ (et on identifie $v$ \`a $\vp$). Si $v$ est une place archim\'edienne, on pose $N(v) = 1$. Souvent dans le texte, on consid\'erera
$$\Sigma_S = \sum_{v \in S} \log N(v).$$
Si $L$ est une extension finie de $K$, et $S'$ est l'ensemble de places de $L$ qui \'etendent celles de $S$, alors
$\Sigma_{S'} = \sum_{\vp \in S} \sum_{\vq | \vp} \log N(\vq) \leq \sum_{\vp \in S} \sum_{\vq | \vp} \frac{f_{\vq}}{f_{\vp}} \log N(\vp) \leq [L:K] \,\Sigma_{S}$, i.e.
\begin{equation}
\Sigma_{S'} \leq [L:K] \Sigma_{S}.  \label{Sigma_S'.Sigma_S}
\end{equation}

On note $P(S) = \{p \in \mathcal{P} / \esp \exists \vp \in S, \vp | p \}$ l'ensemble des caract\'eristiques r\'esiduelles de $S$, et $P$ leur maximum. On montre facilement que
\begin{equation}
\left([K:\rat]\, \card(P(S))\right)^{-1} \Sigma_S \leq \log P \leq \Sigma_{S}, \label{max.car.red}
\end{equation}
 et aussi que 
\begin{equation}
\card(S) \leq [K:\rat] \, \card(P(S)).     \label{cardS}
\end{equation}

\smallskip

%
%

On note $h : \espproj^n(\Qbarre) \longrightarrow [0,\infty)$ la hauteur logarithmique absolue de Weil, et, pour $x$ dans $\espproj^n(K)$, on note $h_K(x) = [K:\rat ]\, h(x)$ la hauteur relative \`a $K$. Avec ces notations,  
\begin{equation}
h_L(x) = [L:K]\, h_K(x).  \label{h_K-h_L}
\end{equation}

Si $\alpha$ est un nombre alg\'ebrique, $h(\alpha)= h(\alpha:1)$.
Si $a,b, c$ sont des \'el\'ements non nuls de $K$ tels que $a+b=c$, on d\'emontre localement que
\begin{equation}
h_{K}(a:b:c) \leq h_{K}(a:c) + [K:\rat]\, \log 2. \label{h(abc).h(ac)}
\end{equation}

\smallskip
%

Pour un point $P=(x_{0}:x_{1}:x_{2})$ de $\espproj^2(K)$, on d\'efinit (cf. \cite{elkies}) le \textit{radical} de $P$ par :
$$\rad_K(P) = \sum_{{\vp }\in I} \log N_{K/\rat}(\vp ),$$
o\`u $I= \{\vp  \in \PK / \espa \card\{ v_{\vp }(x_{0}), v_{\vp }(x_{1}),  v_{\vp }(x_{2})\} \geq 2\}.$ Il d\'epend du corps de rationalit\'e $K$.
Si $a, b$ et $c$ sont des entiers rationnels premiers entre eux,  on a 
$\rad_{\rat}(a:b:c) = \log \prod _{p|abc}  p$.

 Pour tout $r$ appartenant \`a $K$, on pose $P_{r}=(r : 1-r: 1) \in \espproj^2 (K)$ et on peut montrer facilement que
 \begin{equation}
 \rad_K(P_{r}) = \sum_{{\vp }\in H_r} \log N(\vp ), \label{rad(a+b=c)}
 \end{equation}
 o\`u $H_{r} = \{\vp  \in \PK/ \espa v_{\vp }(r)<0 \esp \textrm{ou} \esp v_{\vp }(r)>0 \esp \textrm{ou} \esp v_{\vp }(1-r)>0  \}.$

\bigskip



Le lemme \ref{Sigma_S} sert \`a d\'emontrer le th\'eor\`eme \ref{abc.exp} et le lemme \ref{discriminant}. Ce dernier intervient dans les d\'emonstrations du th\'eor\`eme \ref{SE-abc} et de la proposition \ref{revet.etale}.

\begin{lemme}\label{Sigma_S}

Il existe un nombre r\'eel strictement positif, not\'e dans toute la suite $c_0$, tel que pour tout corps de nombres $K$, tout  ensemble fini $S$ de places ultram\'etriques de $K$, si  $\card(P(S)) \geq 3$, alors 
$$\card(P(S)) \leq c_0\, \frac{\Sigma_{P(S)}}{\log\Sigma_{P(S)} } \leq c_0\, \frac{\Sigma_S}{\log\Sigma_S }.$$

\end{lemme}

\noindent\textbf{D\'emonstration du lemme \ref{Sigma_S}.}

On note $r= \card(P(S))$ et on remarque que $\Sigma_{P(S)}$ est sup\'erieur ou \'egal \`a la somme $T_r$ des logarithmes des $r$ plus petits nombres premiers, et que $\Sigma_S = \sum_{\vp \in S}\log N(\vp) = \sum_{p \in P(S)}\sum_{\vp  |p} f_{\vp } \log p  \geq \sum_{p \in P(S)} \log p = \Sigma_{P(S)}$.

De plus, pour $x \geq e$, la fonction $x \mapsto \frac{x}{\log x}$ est croissante, et comme, pour $r\geq 3$, on a $e \leq T_r \leq \Sigma_{P(S)} \leq \Sigma_S$, alors, d'apr\`es le th\'eor\`eme des nombres premiers, on obtient 
$$\frac{r}{c_0} \leq \frac{T_r}{\log T_r} \leq \frac{\Sigma_{P(S)}}{\log \Sigma_{P(S)}} \leq \frac{\Sigma_S}{\log \Sigma_S}.$$
\hfill $\Box$

\begin{lemme}\label{discriminant}

Soit $L/K$ une extension de corps de nombres. Notons $\delta_{L/K}$ son discriminant relatif et $R=R_{L/K}$ l'ensemble des places de $M_K$ correspondant aux id\'eaux premiers au-dessus desquels l'extension se ramifie, c'est-\`a-dire tels que $v_{\vp }(\delta_{L/K}) > 0$. On d\'esigne par $P(R)$ l'ensemble des caract\'eristiques r\'esiduelles de $R$. Si $\card(P(R)) \geq 3$, alors
$$\log D_{L} \leq [L:K] \left( \log D_{K} + \Sigma_{R} + c_0\, [K:\rat] \log[L:K] \frac{\Sigma_{R}}{\log\Sigma_{R} } \right).$$

\end{lemme}

\noindent\textbf{D\'emonstration du lemme \ref{discriminant}.}

On a $\log D_{L} = \log N_{K/\rat}(\delta_{L/K}) + [L:K] \log D_{K}$ et  $N_{K/\rat}(\delta_{L/K}) = \prod_{\vp  \in R} N(\vp ) ^{\ord_{\vp }(\delta)} $.  D'apr\`es le corollaire \`a la proposition 2, §1 de \cite{serre-chebotarev}, 
$$\ord_{\vp }(\delta_{L/K}) \leq [L:K] - 1 + [L:K]e_{\vp }\frac{\log[L:K]}{\log p},$$
d'o\`u
\begin{displaymath}
\begin{array}{rcl}
\log N_{K/\rat}(\delta_{L/K}) & \leq & \sum_{\vp  \in R} \left( [L:K] + [L:K]e_{\vp } \frac{\log[L:K]}{\log p} \right) \log N(\vp ) \\
               & =    & [L:K] \left( \sum_{\vp  \in R} \log N(\vp ) + \log[L:K] \sum_{\vp  \in R} e_{\vp } f_{\vp } \right) \\
               & =    & [L:K] \left( \Sigma_R + \log[L:K] \sum_{p \in P(R)} \sum_{\vp  |p} e_{\vp } f_{\vp } \right) \\
               & =    & [L:K] \left( \Sigma_R + [K:\rat] \log[L:K]\, \card(P(R)) \right).
\end{array}
\end{displaymath}

On conclut en appliquant le lemme \ref{Sigma_S} \`a l'ensemble des places $R$. \hfill $\Box$

\bigskip

Dans toute la suite, $U$ d\'esigne une courbe alg\'ebrique affine de genre $g$ d\'efinie sur  $K$, $\tilde{U}$ la normalis\'ee de $U$ et $C$ la courbe compl\`ete la contenant. On notera $U_{\infty}= C \setminus \tilde{U} = \{P_{1}, \ldots, P_{t}\}$ l'ensemble des ``points \`a l'infini".

  On notera $h_{U}$ une fonction hauteur d\'efinie sur la courbe $U$ relative \`a un diviseur de degr\'e 1. Par exemple, on peut prendre 
 $$h_{U}= h_{f} =  \frac{1}{\deg(f)}h \circ f,$$
 o\`u $f \in K(U)$ est une fonction non constante sur $U$ et $h$ est la hauteur sur $\espproj^1(\Qbarre)$ d\'efinie pr\'ec\'edemment.
Si $g$ est une autre fonction non constante sur $U$, on a (cf., par exemple, le th\'eor\`eme B.5.9 de \cite{hindry-silverman}), pour tout $x$ de $U(K)$,
\begin{equation}
\left | h_{g}(x) -  h_{f}(x) \right | \leq c \sqrt{h_{f}(x)},  \label{h_f.et.h_g}
\end{equation}
o\`u la constante $c$ ne d\'epend que de la courbe $U$ et des fonctions $g$ et $f$. %

\smallskip


Soit $D \in \Div (C)$ un diviseur tr\`es ample dont le support $\supp (D)$ est l'ensemble des points $U_{\infty}$. 
On notera $U(O_{K,S})$ l'ensemble des points $K$-rationnels de la courbe $U$ qui ont des coordonn\'ees dans l'anneau $O_{K,S}$ des $S$-entiers par rapport au plongement affine $\Phi_{D} : U \hookrightarrow \mathbf{A}^{n}$, donn\'e par $D$.

\smallskip

%
%

Notons $\mathcal{C}$ un mod\`ele de $C$ au-dessus de $\mathcal{S}= \Spec(O_{K,S})$, c'est-\`a-dire que $\mathcal{C} \vers\mathcal{S}$ est un sch\'ema projectif et plat,  tel que sa fibre g\'en\'erique est isomorphe \`a $C$. Alors 
$\mathcal{U} = \mathcal{C} \setminus \overline{\supp(D)}^{Zar} \longrightarrow \mathcal{S}$ est un sch\'ema dont la fibre g\'en\'erique est isomorphe \`a $U$. 

Un point $x$ de $U(K)$ est $S$-entier si et seulement s'il se prolonge en une section $\sigma_x$ au-dessus de $\mathcal{S}$. Dans ce cas l\`a, pour tout id\'eal premier $\vp$ correspondant \`a une place $v$ de $S$, le point $x$ est $\vp $-entier (i.e. $v_{\vp }(x) \geq 0$) et $\sigma_x(\vp )$ appartient \`a la fibre sp\'eciale $\mathcal{U}_{\vp }$.
Nous identifions l'ensemble des points $S$-entiers de la courbe $U$ avec l'ensemble $\mathcal{U}(\mathcal{S})$ des morphismes de $\mathcal{S}$ dans $\mathcal{U}$ au-dessus de $\Spec(K)$, via l'isomorphisme
\begin{displaymath}
\begin{array}{rcl}
U(O_{K,S}) & \simeq & \mathcal{U}(\mathcal{S})\\
      x    & \mapsto & \sigma_x.
\end{array} 
\end{displaymath}

\bigskip


Pour les d\'emonstrations du th\'eor\`eme \ref{SE-abc} et de la proposition \ref{revet.etale}, nous aurons besoin des lemmes \ref{bonne.red.morph} et \ref{ch-weil}. Nous d\'emontrons  ce dernier en utilisant le lemme suivant.

\begin{lemme}\label{revet O_{K,S}}

Soient $K$ un corps de nombres et $S$ un ensemble fini de places de $K$.
Les rev\^etements \'etales de $\Spec (O_{K,S})$ sont des r\'eunions finies de $\Spec (O_{L,S'})$ o\`u $L/K$ est une extension finie, de degr\'e $[L:K]$ inf\'erieur ou \'egal au degr\'e du rev\^etement, non ramifi\'ee en dehors de l'ensemble $S$, et $S' =\{w\in M_L / \esp \exists v \in S, \esp w|v\}$. 

\end{lemme}

\noindent\textbf{D\'emonstration du lemme \ref{revet O_{K,S}}.}

On pose $A= O_{K,S}$ et $Y = \Spec(A)$. Soit $X$ un rev\^etement \'etale de $Y$ et $X = \cup_{i \in I} X_i$ sa d\'ecomposition en composantes connexes. L'ensemble $I$ est fini et, pour tout $i \in I$, $X_i$ est un rev\^etement \'etale connexe de $Y$.

Soit $i \in I$.
On note $L_i = R(X_i)$ l'anneau des fonctions rationnelles de $X_i$ et $Y_i'$ le normalis\'e de $Y$ dans $L_i$. D'apr\`es le corollaire 10.2 de \cite{sga} I, $X_i$ est isomorphe \`a $Y_i'$.
Or, d'apr\`es le corollaire 10.3 de \cite{sga} I, le foncteur $X \mapsto R(X)$ \'etablit une \'equivalence entre la cat\'egorie des rev\^etements \'etales connexes de $Y$ et celle des extensions finies $L/K$ non ramifi\'ees sur $Y = \Spec(O_K) \setminus S$, donc $L_i/K$ est une extension finie non ramifi\'ee en dehors de $S$. 

Comme le foncteur inverse est le foncteur normalisation, alors $X_i$ a pour anneau le normalis\'e $A_i'$ de $A$ dans $L_i$, i.e. $X_i = \Spec(A_i')$ o\`u $A_i' = O_{L_i,S_i'}$, et $S_i'$ est l'ensemble de places du corps $L_i$ qui sont au-dessus de celles de $S$. 

On a donc, $X= \cup_{i \in I} \Spec(O_{L_i, S_i'})$, ce qui d\'emontre le lemme. \hfill $\Box$ 

\medskip

%
%


\begin{lemme}\label{bonne.red.morph}

Soit $f : X \longrightarrow Y$ un morphisme de courbes projectives, d\'efini sur un corps de nombres $K$. 
Il existe un ensemble fini $S_{0}$ de places de $K$ et des mod\`eles $\mathcal{X}$ et $\mathcal{Y}$ de $X$ et de $Y$ au-dessus de $\mathcal{B}= \Spec(O_{K,S_{0}})$ tels que $f$ se prolonge en un morphisme $\tilde{f} : \mathcal{X} \longrightarrow \mathcal{Y}$ et, de plus,

i) si $f$ est fini, il en est de m\^eme de $\tilde{f}$, et dans ce cas, $\deg(f) = \deg(\tilde{f})$,

ii) si, en plus, $f$ est non ramifi\'e en dehors des points $\{P_1, \ldots , P_r\}$ de $Y$, alors, pour tout $\vp  \in \mathcal{B}$, le morphisme $\tilde{f}$ restreint \`a la fibre $\mathcal{X}_{\vp }$ est non ramifi\'e en dehors de $\{ \sigma_{P_1}(\vp ), \ldots, \sigma_{P_r}(\vp ) \}$. Dans ce cas, $\mathcal{X} \setminus \overline{f^{-1}(\{P_1, \ldots , P_r\})}^{Zar} \vers \mathcal{Y} \setminus \overline{\{P_1, \ldots , P_r\}}^{Zar}$ est un rev\^etement \'etale.

\end{lemme}
  
La d\'emonstration du lemme \ref{bonne.red.morph} repose sur le fait que si un nombre, un discriminant, par exemple, est non nul, il en sera de m\^eme pour sa r\'eduction modulo $p$, pour tout nombre premier $p$ sauf un nombre fini. On dit souvent que $X$, $Y$ et $f$ ont  ``bonne r\'eduction" en dehors de $S_{0}$.


\medskip

Le lemme suivant est la version affine du th\'eor\`eme de Chevalley-Weil (cf. \cite{serre} §4.2).

\begin{lemme}\label{ch-weil}

Soient $K$ un corps de nombres, $S_{0}$ un ensemble fini de places de $K$ et $\phi : \mathcal{V} \vers \mathcal{W}$ un rev\^etement \'etale de courbes affines d\'efini sur $\Spec(O_{K,S_{0}})$.  Soient $S_{1}$ un ensemble fini de places de $K$, $y$ un point $S_{1}$-entier de la fibre g\'en\'erique de $\mathcal{W}$ et $x$ un point de $\mathcal{V}_{\eta}(\Qbarre)$ qui rel\`eve $y$.

 Alors le corps de d\'efinition $L=K(x)$ du point $x$ est une extension finie de $K$ de degr\'e $[L:K] \leq \deg(\phi)$, non ramifi\'ee en dehors de l'ensemble $S= S_{1}\cup S_{0}$. De plus, le point $x$ se prolonge en une section $\sigma_x$ au-dessus de $O_{L,S'}$ o\`u $S'= \{w \in M_{L}/ \espa \exists v \in S, \esp w | v \}$. 

\end{lemme}

\noindent\textbf{D\'emonstration du lemme  \ref{ch-weil}.}

Le point $y$ \'etant $S$-entier, o\`u $S= S_{0} \cup S_{1}$, il nous fournit un morphisme de $\Spec(O_{K,S})$ dans $\mathcal{W}$. 
Consid\'erons le produit fibr\'e de $\mathcal{V}$ par $\Spec(O_{K,S})$ au-dessus de $\mathcal{W}$ : 
\begin{displaymath}
\begin{array}{ccc}
\Spec(O_{K,S}) \times _{\mathcal{W}} \mathcal{V} & \longrightarrow & \mathcal{V} \\
pr_1 \downarrow                                             &                 & \downarrow \phi \\
\Spec(O_{K,S})                                   & \longrightarrow & \mathcal{W}.

\end{array}
\end{displaymath}

Le sch\'ema $\Spec(O_{K,S}) \times_{\mathcal{W}} \mathcal{V}$ est obtenu \`a partir de $\mathcal{V}$ par changement de base; le morphisme $\phi : \mathcal{V} \longrightarrow \mathcal{W}$ \'etant un rev\^etement \'etale, il en est donc de m\^eme de $pr_1 : \Spec(O_{K,S}) \times_{\mathcal{W}} \mathcal{V} \longrightarrow \Spec(O_{K,S}) $ (cf. \cite{sga} IX proposition 1.3). De plus, le degr\'e du morphisme $pr_1$ est \'egal \`a celui de $\phi$.

\smallskip

D'apr\`es le lemme \ref{revet O_{K,S}}, $\Spec(O_{K,S}) \times_\mathcal{W} \mathcal{V} = \bigcup_{i \in I} \Spec(O_{L_i,S_i'})$, o\`u, pour tout $i\in I$, l'extension  $L_i/K$ est finie, de degr\'e $[L_i:K] \leq \deg(pr_1) = \deg(\phi)$,  non ramifi\'ee en dehors de l'ensemble $S$, et $S_i' =\{w \in M_{L_i} / \esp \exists v \in S, \esp w | v \}$.

Par ailleurs, le point $x$ de $\mathcal{V}$ \`a valeurs dans $\Qbarre$ correspond \`a un morphisme $\Spec(\Qbarre) \longrightarrow \mathcal{V}$, et l'inclusion $O_{K,S} \subset \Qbarre$, nous fournit un autre morphisme $\iota : \Spec(\Qbarre) \longrightarrow \Spec(O_{K,S})$.
Par la propri\'et\'e universelle du produit fibr\'e, on en d\'eduit qu'il existe un morphisme $\Spec(\Qbarre) \longrightarrow \bigcup_{i \in I} \Spec(O_{L_i,S_i'})$ tel que le diagramme suivant soit commutatif :

\begin{displaymath}
\xymatrix@C=1cm@R=.5cm{
\Spec(\Qbarre) \ar[ddr]_{\iota} \ar@{.>}[dr] \ar[drr]^{x} \\
& \bigcup_i \Spec(O_{L_i,S_i'}) \ar[d] \ar[r] & \mathcal{V} \ar[d]_{\phi} \\
& \Spec(O_{K,S}) \ar[r]                    & \mathcal{W}.                   }
\end{displaymath}

Le spectre d'un corps \'etant r\'eduit \`a un point, le sch\'ema $\Spec(\Qbarre)= \{(0)\}$ s'envoie sur une des composantes connexes de $\bigcup_i \Spec(O_{L_i,S_i'})$, disons sur $\Spec(O_{L,S'})$.

Ainsi, le morphisme $x : \Spec(\Qbarre) \longrightarrow \mathcal{V}$ se factorise par $\Spec(O_{L,S'})$, ce qui revient \`a dire que $x$ se prolonge en une section de $\mathcal{V}$ au-dessus de  $\Spec(O_{L,S'})$. 
Ceci d\'emontre le lemme \ref{ch-weil}. \hfill $\Box$

\bigskip


\section{De ``Siegel uniforme" \`a $abc$.}\label{siegel-abc}
 
 Dans le but d'\'etendre le r\'esultat de L. Moret-Bailly \cite{moret-bailly} aux points $S$-entiers, nous avons cherch\'e \`a formuler un raffinement du th\'eor\`eme  de Siegel.
En suivant les lignes sugg\'er\'ees par N. Elkies (cf. \cite{elkies} et aussi \cite{franken}), nous montrons dans \cite{parma} que la conjecture \ref{abc} impliquerait une version effective du th\'eor\`eme de Siegel, i.e. une borne pour la hauteur des points $S$-entiers, o\`u la d\'ependance en l'ensemble de places $S$ est mise en \'evidence. De plus, si le nombre $c_{\varepsilon,K}$ de la conjecture \ref{abc} \'etait explicite, cette borne serait elle aussi explicite.

On peut essayer de pr\'eciser qu'elle pourrait \^etre la d\'ependance en $K$ de $c_{\varepsilon,K}$. D.W. Masser \cite{masser}, en d\'efinissant un radical {\textit {ramifi\'e}},  montre qu'il faut faire apparaître le logarithme de la valeur absolue $D_{K}$ du discriminant avec un coefficient $>1$. Consid\'erons ce raffinement de la conjecture $abc$, ce qui donnerait, avec la d\'efinition de radical de la section \ref{notations} (cf. \cite{masser} (1.12)), $c_{\varepsilon,K} \leq  2 (1+\varepsilon) \log D_{K} + k_{\varepsilon, [K:\rat]}$. En rajoutant ceci dans la preuve du th\'eor\`eme 1 de \cite{parma}, nous obtenons l'\'enonc\'e suivant.

\begin{thm}\label{th.abc-siegel}

Soient $K$ un corps de nombres, $U$ une courbe alg\'ebrique affine d\'efinie sur $K$ telle que $\chi(U)$ soit strictement n\'egative, $S$ un ensemble fini de places de $K$ et $h_{U,K}$ une fonction hauteur d\'efinie sur $U(K)$ associ\'ee \`a un diviseur de degr\'e 1. 

Supposons vraie la conjecture \ref{abc}, avec  $c_{\varepsilon,K} \leq 2 (1+\varepsilon) \log D_{K} + k_{\varepsilon, [K:\rat]}$, o\`u le nombre $k_{\varepsilon, [K:\rat]}$ ne d\'epend que de $\varepsilon$ et du degr\'e $[K:\rat]$. 

Pour tout $\varepsilon \in ]0,-\chi(U)[$ et tout point $S$-entier $x$ de $U$, il existe une  constante $\kappa$ ne d\'ependant que du degr\'e $[K:\rat]$, de la courbe $U$, du choix de la hauteur $h_{U,K}$ et de $\varepsilon$, telle que 
\begin{equation}
h_{U,K}(x) \leq -\frac{1}{\varepsilon + \chi(U)} \sum_{\vp  \in S} \log N_{K/\rat}(\vp )  -\frac{2 }{\varepsilon + \chi(U)} \log D_K +  \kappa. \label{borne.abc.ram-siegel}
\end{equation}

 De plus, si le nombre $k_{\varepsilon, [K:\rat]}$ \'etait effectif,  il en serait de m\^eme pour  $\kappa$. 
 
 \end{thm}

En nous inspirant d'une part par l'hypoth\`ese du type ``Mordell effectif" fomul\'ee par L. Moret-Bailly dans \cite{moret-bailly} et, d'autre part, par la borne (\ref{borne.abc.ram-siegel}), nous formulons l'hypoth\`ese \ref{se(U,K)}. 
Nous montrons ensuite (th\'eor\`eme \ref{SE-abc}) comment une telle borne, valable sur une courbe donn\'ee, impliquerait une version faible de la conjecture \ref{abc}.

\begin{remarque}\normalfont   

1) L'hypoth\`ese \ref{se(U,K)} est ind\'ependante de la fonction hauteur choisie. (Cf. l'in\'egalit\'e (\ref{h_f.et.h_g})  de la section \ref{notations}.)

2) L'hypoth\`ese \ref{se(U,K)} est trivialement fausse si on ne suppose pas $\chi(U) < 0$. 
\end{remarque}

\smallskip

L'hypoth\`ese \ref{se(U,K)} peut \^etre vue comme l'analogue de plusieurs r\'esultats concernant les corps de fonctions. On y voit apparaître, en particulier, la d\'ependance en l'ensemble de places $S$. (Cf., par exemple, l'appendice de \cite{mathese} pour des \'enonc\'es plus d\'etaill\'es.)

\begin{thm}\label{cdfonctions}

Soient $K=k(B)$ un corps de fonctions, $C$ une courbe  d\'efinie sur $K$ et $h$ une fonction hauteur d\'efinie sur $C$ associ\'ee \`a un diviseur de $C$ de degr\'e 1.
Soit $U$ un ouvert affine de $C$ tel que $\chi(U) <0$. Alors pour toute extension finie $L= k(B')$ de $K$ et tout ensemble fini $S$ de places de $L$, on a, pour tout point $S$-entier $P$ de $U(O_{L,S})$, 
$$h_L(P) \leq c_1 \,\card(S)  + c_2\, (2\,g_L - 2 ) + c_3,$$
o\`u $h_L = [L:K]\, h$ et $g_L = g(B')$ d\'esigne le genre de l'extension $L$ et les nombres $c_1$ et $c_2$ ne d\'ependent que de la courbe $U$ et $c_3$ d\'epend de  $U$, du choix de la hauteur $h$ et du degr\'e $[L:K]$. 

\end{thm}

Remarquons que $2 g_{L}-2$ (respectivement $\card(S)$) est l'analogue pour les corps de fonctions de $\log D_{L}$ (resp. $\sum_{\vp \in S} \log N(\vp)$).
 
Le th\'eor\`eme \ref{cdfonctions} r\'esulte de trois types de r\'esultats connus.
Dans le cas o\`u le genre est nul et que la courbe a au moins trois points \`a l'infini, on peut prendre $c_{1}= c_{2} = 1= \frac{1}{-\chi(U)}$. (C'est essentiellement le th\'eor\`eme de R.C. Mason (lemme 2 de \cite{mason.dio.equa}), analogue de la conjecture $abc$.)  Dans le cas o\`u $g$ et $t$ sont \'egaux \`a 1, on peut prendre $c_{1}=c_{2}= 2 = \frac{2}{-\chi(U)}$. 
(Voir, par exemple, le r\'esultat de M. Hindry et J. Silverman \cite{hindry-silverman.inv}.) Dans le cas o\`u le genre $g$ est sup\'erieur ou \'egal \`a 2 et $t$ est nul, on peut prendre $c_1 =0$ et $c_{2}= \frac{2+\varepsilon}{-\chi(U)}$(cf. les travaux de P. Vojta \cite{vojta}, L. Szpiro \cite{szpiro} et H. Esnault et E. Viehweg \cite{esnault-viehweg}).

D'apr\`es la borne (\ref{borne.abc.ram-siegel}) ainsi que les r\'esultats sur les corps de fonctions ici mentionn\'es, une version optimiste de l'hypoth\`ese \ref{se(U,K)} serait de consid\'erer des constantes $k_1$ et $k_2$ ind\'ependantes de $\delta$, peut-\^etre d\'ependant que de la caract\'eristique d'Euler-Poincar\'e de la courbe. (Le nombre $c_{3}$ d\'ependra de la courbe, du choix de la fonction hauteur et du degr\'e de l'extension de corps de nombres.)

\bigskip

\noindent{\bf D\'emonstration du th\'eor\`eme \ref{SE-abc}.}

Soit $U$ une courbe affine d\'efinie sur un corps de nombres $K$, de genre $g$ et ayant $t$ points \`a l'infini, not\'es $\{P_1, \ldots , P_t \}$, et telle que $\chi(U) <0$. Soit $C$ la courbe projective la contenant. 

\smallskip

On applique le th\'eor\`eme \ref{th.belyi} \`a la courbe $C$. Soit $f : C \longrightarrow \espproj^1$ une fonction de \Belyi \espa de degr\'e $d$, non ramifi\'ee en dehors de $0, 1$ et l'infini, et telle que $f(\{P_1, \ldots , P_t \}) \subset \{0,1, \infty\}$. 
On note $X_f = C \setminus f^{-1}( \{0,1, \infty\})$ la courbe affine obtenue \`a partir de $C$ en lui enlevant les points de ramification de la fonction $f$ et $Y = \espproj^1 \setminus \{0,1, \infty\}$.
D'apr\`es le lemme \ref{bonne.red.morph}, il existe un ensemble fini $S_{0}$ de places de $K$ et des mod\`eles $\mathcal{X}$ et $\mathcal{Y}$ de $C$ et, respectivement de $\espproj^{1}$ au-dessus de $\Spec(O_{K,S_{0}})$ tels que, si on note $\mathcal{X'} = \mathcal{X} \setminus \overline{f^{-1}(\{0, 1, \infty \})}$ et $\mathcal{Y'} = \mathcal{Y} \setminus \overline{\{0, 1, \infty \}} $, alors $\tilde{f} : \mathcal{X'} \vers \mathcal{Y'}$ est un rev\^etement \'etale de degr\'e $d$.

\smallskip

Soient $a,b$ et $c$ des nombres alg\'ebriques non nuls appartenant \`a $K$ tels que $a+b=c$. On pose
$$S_1 = \{\vp  \in \PK / \esp v_{\vp }(a/c)<0 \esp \textrm{ou} \esp v_{\vp }(a/c)>0 \esp  \textrm{ou} \esp  v_{\vp }(b/c)>0 \},$$
de fa\c con \`a ce que $\sum_{\vp \in S_{1}} \log N(\vp) = \rad_{K}(a:b:c)$. 
 
\smallskip

Plongeons $Y = \espproj^1\setminus \{0,1,\infty\}$ dans $\mathbf{A}^3$ \`a l'aide des fonctions $x=T, y=\frac{1}{T}$ et $z=\frac{1}{1-T}$. 
Le point $y=(\frac{a}{c}:1)$ de la courbe $Y$ a pour coordonn\'ees afines $(\frac{a}{c},\frac{c}{a},\frac{c}{b})$.  Comme les nombres $\frac{a}{c},\frac{b}{c}$ et $\frac{c}{a}$ sont des $S_1$-unit\'es, $y$ est un point $S_1$-entier de $Y$.

\smallskip

Soit $x$ dans $X_{f}(\Qbarre)$ un ant\'ec\'edent par $f$ de $y$. 
D'apr\`es le lemme \ref{ch-weil} appliqu\'e \`a $\tilde{f} : \mathcal{X'} \longrightarrow \mathcal{Y'}$,  le corps de d\'efinition $L=K(x)$ du point $x$ est une extension finie de $K$ de degr\'e $[L:K] \leq \deg(f)= d$, non ramifi\'ee en dehors de l'ensemble $S = S_{0} \cup S_{1}$, et, de plus, le point $x$ se prolonge en une section $\sigma_x$ au-dessus de $O_{L,S'}$ o\`u $S'= \{\vq \in M_{L}/ \espa \exists \vp  \in S, \esp \vq | \vp  \}$.

\smallskip

Appliquons l'hypoth\`ese \ref{se(U,K)} \`a la courbe $U$ d\'efinie sur $K$,  \`a $\delta= d$, \`a l'extension de corps de nombres $L/K$ de degr\'e $[L:K] \leq \delta$, \`a la fonction hauteur $h_{U,L} = h_{f,L} = \frac{1}{\deg(f)} (h_L\circ f)$ relative \`a la fonction de \Belyi \espa $f$ et au corps $L$, \`a l'ensemble de places $S'\subset M_L$ et au point $x$ de $X_f \subset U$ qui rel\`eve $y$. On a
\begin{equation}
h_{f,L}(x) \leq k_1\sum_{\vq \in S'}  \log N_{L/\rat}(\vq) + k_2 \log D_L + k_3, \label{appli se(U,h)}
\end{equation} 
o\`u les constantes $k_i(U, h_f, d)_{i\in\{1,2,3\}}$ ne d\'ependent pas du point $x$ ni de $a, b$ et $c$.

\smallskip

Regardons les \'el\'ements qui apparaissent dans cette in\'egalit\'e.

D'apr\`es l'\'egalit\'e (\ref{h_K-h_L}), et parce que $f(x) = (\frac{a}{c}:1)$,
$$h_{f,L}(x) = \frac{1}{\deg(f)} h_L(f(x)) = \frac{[L:K]}{d} h_K(a:c).$$
D'apr\`es l'in\'egalit\'e (\ref{h(abc).h(ac)}), $h_K(a:c) \geq  h_K(a:b:c) - [K:\rat]\log 2$, et comme, $[L:K]\leq d$, alors
\begin{equation}
h_{f,L}(x) \geq \frac{[L:K]}{d} h_K(a:b:c) - [K:\rat]\log 2. \label{haut}
\end{equation}

Gr\^ace \`a la formule (\ref{Sigma_S'.Sigma_S}) on a
\begin{equation}
\Sigma_{S'} \leq [L:K] \,\Sigma_S. \label{rad}
\end{equation}

Pour majorer la valeur absolue du discriminant $D_{L}$ du corps $L$, par rapport \`a celle du corps de base $K$, quitte \`a \'elargir $S_{0}$ pour que $\card(P(S)) \geq 3$, on applique le lemme \ref{discriminant} \`a l'extension $L/K$, qui d'apr\`es le lemme \ref{ch-weil} a pour ensemble de ramification $R_{L/K} \subset S$.
D'o\`u 
\begin{equation}
\log D_L \leq [L:K] \left( \log D_K + \Sigma_S + c_{0} [K:\rat]\log[L:K] \frac{\Sigma_S}{\log \Sigma_S}\right).   \label{discr}
\end{equation}

En rempla\c cant dans l'in\'egalit\'e (\ref{appli se(U,h)}) les in\'egalit\'es (\ref{haut}), (\ref{rad}) et (\ref{discr}), et en remarquant que 
$$\Sigma_S = \sum_{\vp  \in S_0 \cup S_1} \log N(\vp ) \leq \sum_{\vp  \in S_0} \log N(\vp ) + \sum_{\vp  \in S_1}\log N(\vp )  = \rad_K(a:b:c) + k_0,$$
o\`u $k_0 = \Sigma_{S_0}$ ne d\'epend que de la courbe  $U$ et du choix de la fonction de \Belyi, on obtient l'in\'egalit\'e
$$h_K(a:b:c) \leq \gamma_1\, \rad_K(a:b:c) + \gamma_2 \,\frac{\rad_K(a:b:c) + k_0}{\log (\rad_K(a:b:c) + k_0)} + \gamma_3,$$
$$\textrm{o\`u} \esp \gamma_1 = (k_1 + k_2)\,d, \esp \gamma_2 = c_{0}\, k_2\,d\, \log d\, [K:\rat]$$
$$ \textrm{et} \esp \gamma_3 = d\, k_2  \log D_K + d\left(  [K:\rat] \log 2 + k_0 (k_1 + k_2) + k_3 \right).$$ 
 Les constantes $(\gamma_i)_{1\leq i \leq 3}$ et $k_0$ d\'ependent de la courbe $U$ et du corps de nombres $K$ qui ont \'et\'e fix\'es au d\'ebut de la d\'emonstration, mais pas du point $y= (\frac{a}{c}:1)$. Ceci ach\`eve la d\'emonstration du th\'eor\`eme \ref{SE-abc}.\hfill $\Box$

\bigskip


 
\section{$\espproj^1 \setminus \{0,1, \infty\}$, \'equation aux $S$-unit\'es et conjecture $abc$.}\label{droite.proj}

Dans le cas o\`u notre courbe affine $U$ est $\espproj^1 \setminus \{0, 1,\infty\}$, 
la situation est beaucoup plus simple : on n'a pas besoin de toute la force de la conjecture $abc$ pour borner la hauteur des points $S$-entiers de $\droitemoins$; ni d'une hypoth\`ese de type \g {\it Siegel uniforme}" sur la courbe qui fasse intervenir une extension de corps, pour aboutir \`a la conjecture $abc$.

\begin{thm}\label{abc-S-unites}

Soient $K$ un corps de nombres et $\phi$ une fonction positive croissante d\'efinie sur $\mathbf{R}_+$ et d\'ependant \'eventuellement de $K$. On consid\`ere les assertions suivantes  :

i) pour tout triplet $(a,b,c)$ de nombres alg\'ebriques appartenant au corps de nombres $K$, non nuls et tels que $a+b=c$, on a 
$$h_K(a:b:c) \leq \phi (\rad_K(a:b:c)) + \omega \,[K:\rat] \, \log 2;$$

ii) pour tout ensemble fini $S$ de places de $K$ et tout couple $(u, v)$ de $S$-unit\'es v\'erifiant l'\'equation $u+v=1$, on a  
$$ \max\{h_K(u), h_K(v)\} \leq \phi \left( \Sigma_{S}\right).$$

L'assertion i) avec $\omega = 0$ implique l'assertion ii).

L'assertion ii) implique l'assertion i) avec $\omega = 1$.

\end{thm}

 La conjecture \ref{abc} est obtenue \`a partir de l'assertion i) en prenant $\phi(x) = (1 + \varepsilon)x + c_{\varepsilon, K}$.

\begin{thm}\label{S-unites-siegelP1}
Soient $K$ un corps de nombres, $S$ un ensemble fini de places de $K$ et $B_{K,S}$ un nombre r\'eel positif.  En plongeant $U=\droitemoins$ dans l'espace affine $\mathbf{A}^{3}$, on peut choisir la hauteur $h_{U,K}$ de fa\c con \`a avoir une \'equivalence entre :

i) tout couple $(u, v)$ de $S$-unit\'es tel que $u+v=1$, v\'erifie  
$$ \max\{h_K(u), h_K(v)\} \leq B_{K,S};  \textrm{et}$$

ii)  tout point $S$-entier $P$ de $\espproj^1\setminus \{0,1,\infty\}$ v\'erifie
$$h_{U,K}(P) \leq B_{K,S}.$$

\end{thm}

 L'assertion ii) ci-dessus avec $B_{K,S} = c_1 \Sigma_{S} + c_2$ est le cas particulier de l'hypoth\`ese \ref{se(U,K)}, lorsque $U=\espproj^1\setminus\{0, 1, \infty\}$ et $\delta = 1$.

\smallskip

   Si l'on consid\`ere des points entiers sur $\enteros$, i.e. $K=\rat$ et $S = M_{K}^{\infty}$, les th\'eor\`emes \ref{abc-S-unites} et \ref{S-unites-siegelP1} sont vides car l'ensemble $U(\enteros)$ est lui-m\^eme vide, ainsi que l'ensemble des solutions de l'\'equation $u+v=1$. Le probl\`eme n'est int\'eressant que si l'on fait varier $S$.

\bigskip

\noindent{\bf D\'emonstration du  th\'eor\`eme \ref{abc-S-unites}.}

On fixe $K$ et $\phi$ v\'erifiant les hypoth\`eses du th\'eor\`eme. 

\smallskip


\noindent\textbf{L'assertion i) implique l'assertion ii).}\label{demo.abc-S*}

Soient $S$ un ensemble fini de places de $K$ et $u, v$ des $S$-unit\'es de $K$ v\'erifiant l'\'equation $u+v=1$. On pose $u=\frac{a}{c}$ et $v=\frac{b}{c}$ avec $a$, $b$ et $c$ des \'el\'ements non nuls de $K$. Alors
$$u+v=1 \Leftrightarrow  a+b=c.$$

On a $h_K(a:b:c) = h_K(u:v:1) \geq  \max\{h_K(u), h_K(v) \}$ et, d'apr\`es l'\'egalit\'e (\ref{rad(a+b=c)}),
$$\rad_K(a:b:c) = \rad_K(u:1-u:1) = \sum_{\vp  \in H_{u}}\log N(\vp ),$$  
o\`u $H_{u} = \{\vp  \in \PK / \esp v_{\vp}(u)<0  \esp \textrm{ou} \esp v_{\vp}(u)>0  \esp \textrm{ou} \esp v_{\vp}(1-u)>0 \}.$

Comme $u$ et $1-u$ sont des $S$-unit\'es, si la valuation $v_{\vp}$ n'appartient pas \`a $S$,  alors $v_{\vp}(u)$ et $v_{\vp}(1-u)$ sont nulles, et donc la valuation $v_{\vp}$ n'appartient pas \`a $H_u$, c'est-\`a-dire que $(v_{\vp} \in H_u \Rightarrow v_{\vp} \in S) \, \esp \textrm{i.e.} \, \esp H_u \subset S$.
On en d\'eduit que 
$$\rad_K(a:b:c) \leq \sum_{\vp  \in  S}\log N(\vp ) = \Sigma_S.$$

Par hypoth\`ese, $h_{K}(a:b:c) \leq \phi (\rad_{K}(a:b:c))$ et comme la fonction $\phi$ est croissante, 
$$\max\{h_K(u), h_K(v)\} \leq \phi \left( \Sigma_S  \right).$$
\hfill $\Box$

\smallskip


\noindent\textbf{L'assertion ii) implique l'assertion i).}\label{demo.S*-abc}

Soient $a,b$ et $c$ des nombres alg\'ebriques appartenant \`a $K$, non nuls et v\'erifiant $a+b=c$. On pose $S= \{\vp  \in \PK / \esp v_{\vp }(a/c)\ne 0  \esp \textrm{ou} \esp v_{\vp }(b/c)>0 \} $. Alors 
$$u=\frac{a}{c} \esp \textrm{et} \esp v=\frac{b}{c}$$
sont des $S$-unit\'es qui v\'erifient l'\'equation $u+v=1$. On leur applique l'assertion ii). Comme $a+b=c$, d'apr\`es l'in\'egalit\'e (\ref{h(abc).h(ac)}), et parce que $h_{K}(a:c) = h_{K}(u)$, on a  $h_K(a:b:c) \leq \max \{h_K(u), h_K(v)\} + [K:\rat]\log 2$ et, d'apr\`es le (\ref{rad(a+b=c)}),
$$\Sigma_S =  \rad_K (a:b:c).$$ 
On en d\'eduit que $h_K(a:b:c) \leq \phi \left( \rad_K(a:b:c) \right) + [K:\rat]\log2$.
\hfill $\Box$

\bigskip


\noindent\textbf{D\'emonstration du th\'eor\`eme \ref{S-unites-siegelP1}.}\label{demo.S*-siegelP1}

On plonge $U=\espproj^1\setminus\{0, 1, \infty \}$ dans $\textbf{A}^3$ \`a l'aide des fonctions $x=T, y=\frac{1}{T}$ et $z=\frac{1}{1-T}$. On pose aussi $t=1-x$.  Les fonctions $x$ et $t$ prennent aux points $S$-entiers de $U$ des valeurs dans $O_{K,S}$ et on a les relations
$$xy=1, \esp zt=1 \esp \textrm{et} \esp x+t=1,$$
ce qui veut dire que $x$ et $t$ sont des $S$-unit\'es de somme $1$. 

Ainsi, \`a un point $P$ de $U(O_{K,S})$ on associe, en posant $u=x(P)$ et $v=t(P)$, un couple $(u,v)$ de $S$-unit\'es v\'erifiant $u+v=1$, auxquelles nous pouvons appliquer l'assertion i).

R\'eciproquement, \`a un couple $(u,v)$ de $S$-unit\'es v\'erifiant l'\'equation $u+v=1$, on fait correspondre le point $S$-entier $P$ de $\espproj ^1\setminus\{0, 1, \infty\}$ tel que $x(P)=u$, auquel nous appliquons l'assertion ii). Par sym\'etrie, nous pouvons supposer que $h_K(v) \leq h_K(u)$.

En choisissant $h_{U,K}=h_{K} \circ x$, on a le r\'esultat cherch\'e dans les deux situations. \hfill $\Box$

\bigskip

\begin{remarque} \label{u+v=1.Au+Bv=1} \normalfont
On peut consid\'erer une \'equation aux unit\'es plus g\'en\'erale.

Soient $A$ et $B$ des \'el\'ements de $K$ non nuls. 
\'Etant fix\'es le corps de nombres $K$ et la fonction $\phi$, l'assertion ii) du th\'eor\`eme \ref{abc-S-unites} est \'equivalente \`a :

{\it ii') Pour tout ensemble fini $S$  de places de $K$ et tout couple $(u, v)$ de $S$-unit\'es v\'erifiant l'\'equation $Au + Bv =1$, on a  
$$ \max\{h_K(u), h_K(v)\} \leq \phi \left( \Sigma_S + c_{AB} \right) + c'_{AB},$$
o\`u $c_{AB}= \Sigma_{T_{AB}}$,  $T_{AB}= \{\vp  \in \mathcal{P}_{K}/ \espa v_{\vp }(A)\ne 0 \esp \textrm{ou} \esp v_{\vp }(B) \ne 0\}$ 
et $c'_{AB}= h_{K}(A^{-1}:B^{-1}:1)$.}   
\end{remarque}

\noindent\textbf{D\'emonstration de la remarque \ref{u+v=1.Au+Bv=1}.}

Dans le sens r\'eciproque c'est trivial, il suffit de prendre $A=B=1$. On a alors $c_{AB}= c'_{AB}= 0$.
Pour le sens direct, supposons vraie l'assertion ii). Soient $A,B \in K$ non nuls, $S$ un ensemble fini de places de $K$ et $u,v \in O_{K,S}^*$ v\'erifiant l'\'equation $Au+Bv=1$. Posons $S'=S \cup T_{AB}$. Alors $Au$ et $Bv$ sont des $S'$-unit\'es auxquelles on applique la borne de ii). Comme la fonction $\phi$ est croissante, on obtient
$$h_{K}(Au:Bv:1) \leq \phi\left(  \Sigma_{S'} \right) \leq \phi\left(  \Sigma_{S} + \Sigma_{T_{AB}} \right).$$
On conclut en remarquant que $h_{K}(u:v:1) \leq h_{K}(Au:Bv:1) + h_{K}(A^{-1}:B^{-1}:1)$.
\hfill $\Box$

\bigskip


\section{Morphismes de r\'eduction.}\label{meth.baker}


Depuis les travaux de C.L. Siegel, il est connu que la recherche de solutions ($S$-)enti\`eres sur certaines \'equations diophantiennes, peut se ramener \`a la r\'esolution de l'\'equation $u+v=1$ en  ($S$-)unit\'es. C'est ainsi que le r\'esultat (effectif) de A. Baker sur l'\'equation aux unit\'es (obtenu gr\^ace \`a son c\'el\`ebre th\'eor\`eme sur les formes lin\'eaires de logarithmes) a \'et\'e \'etendu \`a d'autres \'equations. Dans \cite{serre} §8.1, J.-P. Serre donne deux \'enonc\'es permettant de d\'eduire la finitude des points entiers d'une courbe, \`a partir d'une autre, \`a l'aide de certains morphismes. Le principe de r\'eduction est explicite de fa\c con tr\`es g\'eom\'etrique et claire dans \cite{serre} et d\'ej\`a pr\'esent implicitement dans \cite{kubert-lang}. Pour le premier \'enonc\'e, on demande que le morphisme soit fini, et pour le second, qu'il soit un rev\^etement \'etale. Ces deux r\'esultats sont ici rendus \g effectifs" au sens o\`u l'on borne la hauteur des points $S$-entiers en fonction de  $S$, ainsi que  du corps de rationalit\'e de la courbe. Nous avons choisi d'exprimer la d\'ependance en l'ensemble $S$ avec  $\Sigma_S = \sum_{v \in S} \log N(v)$, mais aussi avec son cardinal et le maximum des caract\'eristiques r\'esiduelles, $\max P(S)$, qui apparaissent souvent dans la lit\'erature.  (On peut comparer ces quantit\'es entre elles \`a l'aide des in\'egalit\'es de la section \ref{notations}.)

\smallskip

\'Etant fix\'es un corps de nombres $K$, une courbe affine $U$ d\'efinie sur $K$, une fonction hauteur $h_{U}$ d\'efinie sur $U$, et un entier naturel $n \geq 1$, on consid\`ere la propri\'et\'e \g effective" de finitude suivante.

\begin{propriete}{$\mathbf{(P_{U,K,n})}$}\label{prop.finitude}

Il existe une fonction positive $B_{U}$ en cinq variables, croissante en chaque variable, telle que, pour toute extension $L/K$ de degr\'e $[L: K] \leq n$, tout ensemble fini $T$ de places de $L$ et tout point $x$ de $U(O_{L,T})$, on ait
$$h_{U}(x)\leq B_{U}(\Sigma_{T}, \card(T), \max P(T), \log D_{L}, [L:\rat]).$$
\end{propriete}

Cette propri\'et\'e ne d\'epend pas de la fonction hauteur choisie (cf. l'in\'egalit\'e (\ref{h_f.et.h_g})). En revanche, la fonction $B_{U}$ en d\'epend.


\begin{prop}\label{morph.fini}
Soient $\phi : X \vers Y$ un morphisme fini de courbes alg\'ebriques affines  d\'efini sur un corps de nombres $K$ et $h_{X}$ et $h_{Y}$ des  hauteurs sur $X$ et $Y$. Supposons que la courbe $Y$ v\'erifie la propri\'et\'e  $\mathbf{(P_{Y,K,n})}$ pour $n = 1$ et une fonction $B_{Y}$.

Alors la courbe $X$ v\'erifie la propri\'et\'e  $\mathbf{(P_{X,K,n})}$ pour $n=1$ et une fonction $B_{X}$. De plus, on peut choisir la fonction  $B_{X}$ v\'erifiant :

il existe des r\'eels $u_{\phi}, v_{\phi}$ et $w_{\phi}$  ne d\'ependant que de $\phi : X \longrightarrow Y$, et $\gamma >0$,  tels que  
$$ B_{X}(u, v, w, z, d)= \gamma \, B_{Y}(u + u_{\phi}, v +v_{\phi}, w + w_{\phi}, z, d).$$

\end{prop}


\begin{prop}\label{revet.etale}
Soit $\phi : X \vers Y$ un morphisme de courbes alg\'ebriques affines de degr\'e $d_{\phi}$, d\'efini sur un corps de nombres $K$. Supposons $\phi$ fini, surjectif et non ramifi\'e.
Soient $h_{X}$ et $h_{Y}$ des fonctions hauteurs sur $X$ et $Y$, respectivement. Soit $n \geq d_{\phi}$ et supposons que la courbe $X$ v\'erifie la propri\'et\'e $\mathbf{(P_{X,K,n})}$ pour une fonction $B_{X}$.

Alors  la courbe $Y$ v\'erifie la propri\'et\'e  $\mathbf{(P_{Y,K,\delta})}$ pour  $\delta = \frac{n}{d_{\phi}}$ et une fonction $B_{Y}$. De plus, on peut choisir la fonction $B_{Y}$ v\'erifiant :
 
il existe des  r\'eels $u_{\phi}, v_{\phi}$ et $w_{\phi}$ ne d\'ependant que de $\phi : X \longrightarrow Y$, et $\gamma>0$,  tels que 
$$B_{Y}(u, v, w, z, d) = \gamma\, B_{X}(d_{\phi}(u + u_{\phi}), d_{\phi}(v+ v_{\phi}), w + w_{\phi}, \gamma_{d}, d_{\phi} \,d),$$
o\`u $\gamma_{d} = d_{\phi} (z + u + u_{\phi} + c_{0}\, d \,\log(d_{\phi}) \frac{u + u_{\phi}}{\log(u +  u_{\phi})})$ et $c_{0}$ est la constante absolue du lemme \ref{Sigma_S}.

\end{prop}

Sous les hypoth\`eses de la proposition \ref{revet.etale}, si $Y$ v\'erifie $\mathbf{(P_{Y,K,n})}$, alors $X$ v\'erifie $\mathbf{(P_{X,K,n})}$.

\begin{remarque}\normalfont   Dans les propositions \ref{morph.fini} et \ref{revet.etale}, on peut choisir $u_{\phi} = \Sigma_{S_{\phi}}, v_{\phi} = \card(S_{\phi})$ et $w_{\phi} = \max P(S_{\phi})$, o\`u $S_{\phi}$  est l'ensemble fini des places en lesquelles les courbes $X$ ou $Y$ ou le morphisme $\phi$ ont \g mauvaise r\'eduction". Le nombre $\gamma$ d\'epend du  morphisme $\phi : X \longrightarrow Y$, ainsi que des choix des hauteurs $h_{X}$ et $h_{Y}$.
\end{remarque}

\begin{remarque} \label{phi^{-1}} \normalfont 

Notons toutefois que nous ne pouvons pas appliquer la proposition  \ref{morph.fini} \`a la fonction de \Belyi \esp associ\'ee \`a une courbe $U$ pour borner la hauteur de ses points $S$-entiers \`a partir des bornes inconditionnelles connues pour la hauteur des points $S$-entiers de $\droitemoins$. En effet, le th\'eor\`eme de \Belyi \esp nous permet d'appliquer la proposition \ref{morph.fini} \`a l'ouvert $X = X_{f}$, et non pas \`a $U$, car il nous donne l'inclusion $f(U_{\infty}) \subset \{0,1,\infty\}$, mais pas l'\'egalit\'e.

\end{remarque}

\medskip

\noindent\textbf{D\'emonstration de la proposition  \ref{morph.fini}.}

Soient $h_{X}$ et $h_{Y}$ deux fonctions hauteurs d\'efinies respectivement sur la courbe $X$ et la courbe $Y$, associ\'ees \`a des diviseurs de degr\'e 1. Alors la fonction $\frac{1}{\deg(\phi)} h_{Y} \circ \phi$ est une hauteur sur $X$ associ\'ee \`a un diviseur de degr\'e 1. 
D'apr\`es la th\'eorie des hauteurs (cf. l'in\'egalit\'e (\ref{h_f.et.h_g})), pour tout point rationnel $x$ de $X$, on a 
$$h_{X}(x) \leq \frac{1}{\deg(\phi)} h_{Y}(\phi(x)) + O(\sqrt{ h_{Y}(\phi(x))}),$$
o\`u la constante implicite ne d\'epend que des courbes $X$ et $Y$, des hauteurs $h_{X}$ et $h_{Y}$ et du morphisme $\phi$.

Si $S$ est un ensemble fini de places du corps $K$ et si le point $x$ est $S$-entier, alors $\phi(x)$ est un point $S'$-entier de la courbe $Y$, o\`u $S' = S \cup S_{\phi}$ pour un ensemble fini $S_{\phi}$ ne d\'ependant que du morphisme $\phi : X \longrightarrow Y$. 
En remarquant d'une part que $\Sigma_{S'} \leq \Sigma_{S} + \Sigma_{S_{\phi}}$, que  $\card(S') \leq \card(S) + \card(S_{\phi})$ et que $\max P(S') \leq \max P(S) + \max P(S_{\phi})$, et d'autre part que la fonction $B_{Y}$ est positive et croissante, on obtient
$$h_{Y}(\phi(x)) \leq B_{Y}(\Sigma_{S} + \Sigma_{S_{\phi}}, \card(S) + \card(S_{\phi}), \max P(S) + \max P(S_{\phi}), \log D_{K}, [K:\rat]);$$
ce qui nous permet de conclure. 
\hfill $\Box$


\bigskip

\noindent\textbf{D\'emonstration de la proposition  \ref{revet.etale}.}

Soit $L/K$ une extension de corps de degr\'e $[L:K] \leq \frac{n}{d_{\phi}}$.

Notons $\overline{X}$ (respectivement $\overline{Y}$) la compl\'et\'ee de $X$ (resp. de $Y$) et $X_{\infty} = \overline{X}\setminus X$ (resp. $Y_{\infty}$) les points \g \`a l'infini".  Le morphisme $\phi : X \vers Y$ s'\'etend en un morphisme  (toujours not\'e $\phi$) de $\overline{X}$ vers $\overline{Y}$. Comme il est fini, $\phi^{-1}(Y_{\infty}) = X_{\infty}$. D'apr\`es le lemme \ref{bonne.red.morph}, il existe un ensemble fini $S_{\phi}$ de places de $L$ en dehors duquel $\phi$ s'\'etend en un rev\^etement \'etale $\tilde{\phi} : \mathcal{X} \setminus \overline{X_{\infty}} \vers \mathcal{Y} \setminus \overline{Y_{\infty}}$.

Soient  $S$ un ensemble fini de places de $L$ et $y$ un point $S$-entier de $Y$. Comme le morphisme $\phi$ est surjectif, on peut relever le point $y$ en un point $x$ de $X(\Qbarre)$. 

D'apr\`es le lemme \ref{ch-weil}, le corps de rationalit\'e $M= L(x)$ de $x$ est une extension finie, de degr\'e $[M:L] \leq d_{\phi}$ et non ramifi\'ee en dehors de l'ensemble $S \cup S_{\phi}$. De plus, le point $x$ est $S'$-entier, o\`u $S'$ est l'ensemble de places de $M$ au-dessus de celles de $S \cup S_{\phi}$.

Soit $h_{X}$ une fonction hauteur d\'efinie sur la courbe $X$.  Comme $[M:K] \leq d_{\phi}\, \frac{n}{d_{\phi}} = n$, appliquons l'hypoth\`ese \`a la fonction hauteur $h_{X}$, \`a l'extension $M$ de $K$,  \`a l'ensemble de places $S'$ et au point $x$ qui rel\`eve $y$. On a 
$$h_{X}(x) \leq B_{X}(\Sigma_{S'}, \card(S'), \max P(S'), \log D_{M}, [M:\rat]).$$
En appliquant l'in\'egalit\'e (\ref{Sigma_S'.Sigma_S}), on a $\Sigma_{S'} \leq [M:L] \Sigma_{S \cup S_{\phi}} \leq d_{\phi} (\Sigma_{S} + \Sigma_{\phi})$. De plus, $\card (S') \leq [M:L] \,\card (S \cup S_{\phi}) \leq d_{\phi}\left(\card (S) + \card (S_{\phi}) \right)$ et $\max P(S') = \max P(S\cup S_{\phi}) \leq \max P(S) + \max P(S_{\phi})$. Gr\^ace au lemme \ref{discriminant}, on majore la valeur absolue du discriminant de l'extension $M$ en fonction de celle du corps de base $L$. On a $\log D_{M} \leq \gamma_{L}$, o\`u
$$\gamma_{L} =d_{\phi} \left(\log D_{L} + \Sigma_{S} + \Sigma_{S_{\phi}} + c_{0} \, [L:\rat] \,\log(d_{\phi} ) \frac{\Sigma_{S} +\Sigma_{S_{\phi}}}{\log (\Sigma_{S} + \Sigma_{S_{\phi}})}  \right)$$
et $c_{0}$ est la constante du lemme \ref{Sigma_S}.
D'o\`u
$$h_{X}(x) \leq  B_{X}(d_{\phi}\,(\Sigma_{S} + \Sigma_{S_{\phi}}), d_{\phi}\,(\card (S) + \card (S_{\phi})), \max P(S) + \max P(S_{\phi}),  \gamma_{L}, d_{\phi}\,[L:\rat]). $$

Par ailleurs,  si $h_{Y}$ est une fonction hauteur d\'efinie sur la courbe $Y$ et associ\'ee \`a un diviseur de degr\'e 1, alors $\frac{1}{d_{\phi}}(h_{Y} \circ \phi)$ est une fonction hauteur d\'efinie sur la courbe $X$ associ\'ee \`a un diviseur de degr\'e 1 et on a (cf. l'in\'egalit\'e (\ref{h_f.et.h_g})) 
$$h_{Y}(y) \leq d_{\phi}\,h_{X}(x) + O(\sqrt{h_{X}(x)}),$$
o\`u la constante implicite ne d\'epend que de $X, Y, h_X, h_{Y}$ et $\phi$, et on peut conclure.\hfill $\Box$

 \begin{cor}\label{coro}
Les th\'eor\`emes  \ref{abc-S-unites} et \ref{S-unites-siegelP1} impliquent le th\'eor\`eme \ref{SE-abc}.
 
\end{cor}

\noindent\textbf{D\'emonstration du corollaire \ref{coro}.}
 
 L'hypoth\`ese \ref{se(U,K)} peut \^etre vue comme la collection des propri\'etes $\mathbf{(P_{U,K,\delta})}$ pour tout $\delta \geq 1$, avec la fonction $B_{U}(u, v, w, z, d) = k_{1}\,u + k_{2}\, z + k_{3}$.
  
   D'apr\`es les th\'eor\`emes \ref{abc-S-unites} et \ref{S-unites-siegelP1}, qui montrent ensemble qu'un \'enonc\'e du type {\it \g Siegel Uniforme"} pour $\droitemoins$ impliquerait la conjecture $abc$, le th\'eor\`eme \ref{SE-abc} peut \^etre obtenu en appliquant la proposition \ref{revet.etale} \`a une fonction de \Belyi \esp $f$ associ\'ee \`a la courbe $U$. Si la courbe $U$ est d\'efinie sur un corps $K$, en prenant $\delta = \deg(f)$, on obtient un \'enonc\'e allant dans le sens de la conjecture $abc$ valable sur $K$. (Cf. le th\'eor\`eme \ref{abc.exp}.)
\hfill $\Box$

\bigskip


\section{Applications.}\label{applications}

Commen\c cons par \'enoncer le th\'eor\`eme de Y. Bugeaud et K. Gy\"ory \cite{bugeaud-gyory} (concernant l'\'equation aux unit\'es) et le th\'eor\`eme  de Yu. Bilu \cite{bilu} (version effective du th\'eor\`eme de Siegel pour les rev\^etements galoisiens de la droite projective), qui nous permettent d'obtenir un r\'esultat inconditionnel : le th\'eor\`eme \ref{abc.exp}. 

\begin{thm}{\bf(Bugeaud-Gy\"ory)}\label{thm.bugeaud-gyory}
Soient $K$ un corps de nombres, $H$ un nombre r\'eel $\geq e$, $A$ et $B$ des \'el\'ements non nuls de $K$ tels que $\max \{h(A), h(B) \} 
 \leq \log H$ et $T$ un ensemble fini de places de $K$ contenant les places archim\'ediennes.
Les solutions $u, v$ de l'\'equation $Au+Bv=1$ appartenant \`a $O_{K,T}^*$ v\'erifient
$$\max\{h(u), h(v) \} \leq   \gamma P^{[K:\rat]} R_T\, (\log^+ R_T) (\log^+ (P R_T)/\log^+ P)  \log H,$$
o\`u $\gamma$ est un nombre r\'eel d\'ependant de $[K:\rat]$ et de $\card(T)$, $P$ est le maximum des caract\'eristiques r\'esiduelles de $T$, $R_T$ est le $T$-r\'egulateur et $\log^+ (.) $ est une notation pour $\max\{ \log ., 1\}$.

\end{thm}

Pr\'ecis\'ement, $\gamma = c_K^{\card(T) +1}\, (\card(T) )^{5\, \card(T)+ 10}$, o\`u $c_K \geq 1$ ne d\'epend que de $[K:\rat]$.

\begin{thm}{\bf(Bilu)}\label{thm.bilu}
Soient $K$ un corps de nombres, $C$ une courbe alg\'ebrique projective de genre $g \geq 1$, d\'efinie sur $K$ et $x \in K(C)$ une fonction non constante telle que $x : C \vers \espproj^1$ soit un rev\^etement galoisien (i.e. $\Qbarre(C)/ \Qbarre(x)$ et une extension galoisienne).
Pour tout ensemble fini $S$ de places de $K$ contenant les places archim\'ediennes, on pose 
$$C(x, K, S) = \{P \in C(K)/ \espa x(P) \in O_{K,S} \}.$$

Soient $y \in K(C)$ telle que $K(C) = K(x,y)$ et $f(X,Y) \in K[X,Y]$ un polyn\^ome s\'eparable non nul tel que $f(x,y)=0$. On pose $m= \deg_{X}(f)$ et $n = \deg_{Y}(f)$. 
Pour tout $P$ de $C(x,K,S)$ on a 
$$h(x(P)) \leq P^{N_1 [K:\rat]}\left( D_K \prod_{\vp \in S} N_{K/\rat}(\vp) \right)^{N_2} e^{\psi},$$
o\`u $\psi = 100\,\card(S)N_2\left( \log(N\, \card(S)) + O(1) \right) + [K: \rat] N_3 \left( h(g) + O(N) \right)$, 
$N = \max\{m,n,3\}$, $N_1= \max\{n^{5}, 16 n^{2}m^{2}, 256 m^{3}\}$, $N_2 = \max \{n^{4}, 10 m^{2}n\}$ et $N_3= \max \{ mn^{7}, 500 m^{2}n^{4}\}$, 
  $h(g)$ d\'esigne la hauteur du polyn\^ome $g$, \`a savoir, la hauteur du point de l'espace projectif d\'efini par les coefficients de $g$, et  
$P$ est le maximum des caract\'eristiques r\'esiduelles de l'ensemble $S$.

\end{thm}

\smallskip

 L'in\'egalit\'e du th\'eor\`eme \ref{abc.exp} est du m\^eme ordre (exponentiel) que celle des r\'esultats pr\'ec\'edents de Stewart-Tijdeman \cite{stewart-tijdeman} et Stewart-Yu \cite{stewart-yu}, et les constantes sont de qualit\'e inf\'erieure; elle a n\'eanmoins l'avantage d'\^etre valable pour tout corps de nombres. 


\begin{thm}{\bf(Stewart-Yu)}\label{thm.stewart-yu}
Il existe une constante $\eta>0$ effectivement calculable telle que, pour tout triplet $(a,b,c)$ d'entiers naturels positifs, premiers entre eux et v\'erifiant $a+b=c$, on ait 
$$h(a:b:c) < \eta \espa \rad_{\rat}(a:b:c)^3 \exp \left( \frac{1}{3} \rad_{\rat} (a:b:c) \right).$$
\end{thm}

Le th\'eor\`eme \ref{thm.stewart-yu} correspond \`a l'assertion i) du th\'eor\`eme \ref{abc-S-unites} o\`u le corps de nombres est $\rat$, la fonction $\phi$ est d\'efinie par $\phi(x) = \eta x^3 \exp(\frac{1}{3}x)$ et $\omega = 0$.


\subsection{Vers $abc$.}

Pour obtenir le th\'eor\`eme \ref{abc.exp}, on majore les termes de la borne du th\'eor\`eme \ref{thm.bugeaud-gyory} faisant intervenir l'ensemble $T$, \`a savoir $P$, $R_T$ et $\card(T)$, en fonction de $D_{K}$ et de $\Sigma_T$, et on applique le th\'eor\`eme \ref{abc-S-unites} \`a la fonction $\phi$ de $\Sigma_T$ ainsi obtenue.

\smallskip

\noindent{\bf D\'emonstration du th\'eor\`eme \ref{abc.exp}.}

Soit $T$ un ensemble fini de places du corps $K$ et $u$ et $v$ des $T$-unit\'es v\'erifiant l'\'equation $u+v = 1$, auxquelles on applique le th\'eor\`eme \ref{thm.bugeaud-gyory}.
On d\'esigne par $c_{i}$ des nombres r\'eels ne d\'ependant que du degr\'e $[K:\rat]$ et on pose  $S= T\cap \mathcal{P}_{K}$.

La remarque 1 de \cite{bugeaud-gyory} nous donne $\frac{\log^{+}(PR_{T})}{\log^{+}P} \leq 2\,\log^{+}R_{T}$.
Quitte \`a \'elargir $S$ (ce qui modifie l\'eg\'erement les nombres $c_{i}$ ne d\'ependant ni de $S$ ni de $D_{K}$), nous pouvons supposer que $\card(P(S)) \geq 3$ et $\Sigma_{S}\geq e$. Alors $\log^{+}R_{T} = \log R_{T}$ et $R_{T}\,\log^{+}R_{T} \, \frac{\log^{+}(P R_{T})}{\log^{+} P} \leq 2\,R_{T}^{2}$.
En appliquant succesivement le lemme 3 de \cite{bugeaud-gyory} et le lemme 8 de \cite{bugeaud.super}, on obtient
$$R_{T} \leq R_{K}\,\mathfrak{h}_{K}\, \prod_{\vp \in S}\log N(\vp) \leq c_{1} \sqrt{D_{K}}\,(\log D_{K})^{[K:\rat] -1} \, \prod_{\vp \in S}\log N(\vp),$$
 o\`u $\mathfrak{h}_{K}$ d\'esigne le nombre de classes. On a alors     
\begin{equation}
\max \{h_{K}(u), h_{K}(v)\} \leq c_{2}\,\gamma_{(T, [K:\rat])}\,P^{[K:\rat]} \, D_{K}\,(\log D_{K})^{2[K:\rat] -2}  \,(\prod_{\vp \in S} \log N(\vp) )^{2}. \label{borne.u+v.4var}
\end{equation}
D'apr\'es l'in\'egalit\'e (\ref{max.car.red}) nous avons
$P^{[K:\rat]} \leq \exp \{[K:\rat]\, \Sigma_{S}\}$.
En appliquant l'in\'egalit\'e (\ref{cardS}) et le lemme \ref{Sigma_S} on obtient 
$\card(S)\leq c_{0}\,[K:\rat]\,\frac{\Sigma_{S}}{\log \Sigma_{S}}$ et comme $\card(T) \leq \card(S) + [K:\rat]$, alors $\gamma \leq \exp \{c_{3}\, \Sigma_{S}\}$. D'apr\`es l'in\'egalit\'e arithm\'etico-g\'eom\'etrique, 
$\prod_{\vp \in S} \log N(\vp) \leq (\Sigma_{S})^{\card (S)}$. 
D'o\`u l'existence de r\'eels  $\gamma_{1}$ et $\gamma_{3}$ effectivement calculables et ne d\'ependant que de $[K:\rat]$, tels que
$$\max\{h_{K}(u), h_{K}(v) \} \leq \exp \{ \gamma_{1} \Sigma_{S} + \log D_{K} + (2\,[K:\rat]-2) \log \log D_{K}+ \gamma_{3} \}.$$

On conclut en appliquant le th\'eor\`eme \ref{abc-S-unites} \`a la fonction $\phi(x) = \exp \{ \gamma_{1}\, x + \gamma_{2}\, \log D_{K} + \gamma_{3} \}$, avec  $\gamma_{2} = 2[K:\rat] -1$.
\hfill $\Box$

\begin{remarque}\label{bilu-abc.exp}
Bien que la borne de l'hypoth\`ese \ref{se(U,K)} soit plus forte que celle du th\'eor\`eme \ref{thm.bilu}, nous pouvons d\'eduire de ce dernier, en reprenant la d\'emonstration du th\'eor\`eme \ref{SE-abc}, le th\'eor\`eme \ref{abc.exp}. Nous obtenons ainsi une constante absolue $\gamma_{2}$, ind\'ependante du degr\'e $[K:\rat]$.
\end{remarque}

\noindent{\bf D\'emonstration de la remarque \ref{bilu-abc.exp}.}

Soient $K$ un corps de nombres et $a,b$ et $c$ des \'el\'ements non nuls de $K$ tels que $a+b=c$. 
Posons $S_{1} =\{\vp \in \mathcal{P}_{K} / \espa v_{\vp}(\frac{a}{c}) \ne 0 \esp \textrm{ou} \esp v_{\vp}(\frac{b}{c}) > 0 \}$, de fa\c con \`a ce que $\Sigma_{S_{1}} = \rad_{K}(a:b:c)$ et que $(a:c) \in (\droitemoins) (O_{K,S})$.

Soient $U$ la courbe d'\'equation affine $y^{2}= x^{3} -x$ et $C$ la courbe projective correspondant \`a l'\'equation homog\`ene $Y^{2}Z = X^{3}- XZ^{2}$.  Notons $P_{\infty} = (0:1:0)$ le point \`a l'infini. 

Soit $f$ une fonction de \Belyi \esp associ\'ee \`a la courbe $C$ telle que $f(P_{\infty}) \subset \{0,1,\infty\}$. Posons $d = \deg(f)$. Soit $p$ un point de $C\setminus f^{-1}(\{0,1,\infty\}) \subset U$ qui rel\`eve le point $S$-entier $(a:c)$.  

D'apr\`es les lemmes \ref{bonne.red.morph} et \ref{ch-weil}, il existe un ensemble fini $S_{0}$ d'id\'eaux premiers de $K$, tel que le corps de rationalit\'e $L=K(p)$ du point $p$ soit une extension finie de degr\'e $[L:K] \leq d$, non ramifi\'ee en dehors de l'ensemble $S=S_{0}\cup S_{1}$ et le point $p$ soit $S'$-entier, o\`u $S'$ d\'esigne l'ensemble de places de $L$ qui sont au-dessus de celles de $S$.

Quitte \`a \'elargir l'ensemble $S_{0}$, supposons que $\card(P(S)) \geq 3$.

Au lieu de l'hypoth\`ese {\it \g Siegel Uniforme" $(U,K)$}, appliquons ici le th\'eor\`eme \ref{thm.bilu} \`a la courbe $C$, au corps $L$, \`a la fonction rationnelle $x$ et au point $S'$-entier $p$. (Avec nos notations, l'ensemble $U(O_{L,S'})$ des points $S'$-entiers de $U$ est $\{(x,y) \in (O_{L,S'})^{2} /\esp  y^{2} = x^{3} -x \}$ et cet ensemble est inclus dans celui consid\'er\'e par Yu. Bilu, $C(x, L, S')$. De plus, la hauteur du polyn\^ome d\'efinissant notre courbe est nulle.) Ainsi
$$h_{x}(p) \leq \exp \left\{ N_2 (\Sigma_{S'} + \log D_L ) + N_1[L:\rat] \log P + \psi \right\},$$
o\`u  $\psi = 400\,N_{2}\card(S')\left( \log( \card(S')) + N_{4} \right) + [L:\rat] N_3$, $P$ est le maximum des caract\'eristiques r\'esiduelles de l'ensemble $S'$ et $N_1$, $N_2, N_{3}$ et $N_4$ sont des constantes abolues, effectivement calculables.
 
Ramenons-nous au corps $K$ et \`a l'ensemble $S$.

Comme $[L:K] \leq d$, d'apr\`es l'in\'egalit\'e (\ref{Sigma_S'.Sigma_S}) on a $\Sigma_{S'} \leq d \,\Sigma_S$, et d'apr\`es le lemme \ref{discriminant},
$$\log D_L \leq d \left( \log D_K + \Sigma_S + c_{0}\,[K:\rat] \log d \frac{\Sigma_S}{\log\Sigma_S } \right).$$
Par ailleurs, l'in\'egalit\'e (\ref{max.car.red}) nous majore le maximum $P$ des caract\'eristiques r\'esiduelles de $S'$ (ou de $S$) : $\log P \leq  \Sigma_S$.
L'in\'egalit\'e (\ref{cardS}),  ainsi que le lemme \ref{Sigma_S}, appliqu\'es \`a l'ensemble $S'$, nous permettent de majorer le cardinal de $S'$ : $\card(S')  \leq c_{0}\, d\, [K:\rat] \frac{\Sigma_{S}}{\log \Sigma_{S}}$, ce qui nous permet de majorer $\psi$. On a
$$h_{x}(p) \leq \exp \left\{ n_{1} \Sigma_{S} + n_{2} \frac{\Sigma_{S}}{\log \Sigma_{S}} + n_{3}\, d \log D_{K} + n_{4} \right\},$$
o\`u $n_{1}, n_{2}$ et $n_{4}$ sont des polyn\^omes en $d, \log d$ et $[K:\rat]$ de degr\'e inf\'erieur ou \'egal \`a 1 en chacune des variables.

Ramenons-nous maintenant \`a la hauteur et au radical de $(a:b:c)$.

Remarquons que $\Sigma_S \leq \Sigma_{S_{1}} + \Sigma_{S_{0}} = \rad_K(a:b:c) + \Sigma_{S_0}$ et que $S_0$ ne d\'epend que de la courbe $U$ et de la fonction de \Belyi \esp choisies au d\'ebut.

De plus, on a $h_{f}(p) = \frac{h(f(p))}{d} = \frac{h(a:c)}{d} = \frac{h_{K}(a:c)}{d\,[K:\rat]}$ et, d'apr\`es l'in\'egalit\'e (\ref{h(abc).h(ac)}), 
$$h_{f}(p) \geq \frac{h_{K}(a:b:c)}{d\,[K:\rat]} - \frac{\log 2}{d}.$$
Par ailleurs, d'apr\`es l'in\'egalit\'e (\ref{h_f.et.h_g}), $h_{f}(p) \leq h_{x}(p) + n_{5}\sqrt{h_{x}(p)}$, avec $n_{5}$ d\'ependant des fonctions $x$ et $f$ et de la courbe $U$, donc
$$h_{K}(a:b:c) \leq d\,[K:\rat] \left(h_{x}(p) + n_{5} \sqrt{h_{x}(p)}\right) + [K:\rat] \, \log 2.$$
On obtient

$$h_K(a:b:c) \leq \exp\{\gamma_{1} \,\rad_K(a:b:c) + \gamma_{2} \log D_{K} + \gamma_{3} \},$$
o\`u le nombre r\'eel  $\gamma_{1}$ d\'epend de $d$ et de $[K:\rat]$, $\gamma_{2} = n_{6} \, d$ avec $n_{6}$ une constante absolue et $\gamma_{3}$ d\'epend de $x, f, U, d$ et $[K:\rat]$. Ceci ach\`eve la d\'emonstration. \hfill $\Box$

\medskip

Au lieu de la courbe d'\'equation $y^{2}= x^{3}-x$, on aurait pu prendre n'importe quel rev\^etement galoisien de la droite projective.


\medskip

\subsection{Courbes de genre nul.}\label{genre.nul}

Dans ce paragraphe, $K$ est un corps de nombres, $T$ un ensemble fini de places de $K$ contenant les places archim\'ediennes, $P$ est le maximum des caract\'eristiques r\'esiduelles de $T$, et on pose $S= T\cap \mathcal{P}_{K}$, $d=[K:\rat]$, $t=\card(T)$ et $s=\card(S)$.  

\begin{cor}\label{cor.P1}
 Tout point $T$-entier $x$ de $\droitemoins$ v\'erifie :
$$h_{K} (x) \leq  c_{d}^{t}\,t^{5(\,t + 2)}\,P^{d}  \,(\log P )^{2\,s} \, D_{K}\,(\log D_{K})^{2\,d -2},$$
o\`u  le nombre $c_{d}$ ne d\'epend que du degr\'e $d$.
\end{cor}

\medskip

\noindent{\bf D\'emonstration du corollaire \ref{cor.P1}.}

On applique le th\'eor\`eme \ref{S-unites-siegelP1} \`a la borne (\ref{borne.u+v.4var}), obtenue pour la hauteur des solutions en $T$-unit\'es de l'\'equation $u+v=1$ gr\^ace au th\'eor\`eme \ref{thm.bugeaud-gyory}. Puis on remarque que $\prod_{\vp \in S} \log N(\vp) \leq  \prod_{\vp \in S} f_{\vp} \,\log P \leq (\log P)^{s}\, \prod_{p\in P(S)} (\sum_{\vp | p} f_{\vp})^{\card\{\vp | p\}} \leq (\log P)^{s} \, d^{d\,\card(P(S))} \leq (d^{d} \, \log P)^{s}$. 
\hfill $\Box$

\begin{remarque}\label{g=0}
Si $U$ est une courbe de genre nul ayant au moins trois points \`a l'infini, on peut appliquer la proposition \ref{morph.fini} au morphisme $\phi : U \longrightarrow \droitemoins$ qui identifie $C$ \`a $\espproj^{1}$ et envoi $U_{\infty}$ sur $\{0,1,\infty\}$, et \`a la borne obtenue pour $\droitemoins$ dans le corollaire \ref{cor.P1}, pour obtenir, pour tout point  $x$ de $U(O_{K,T})$,
\begin{equation}
h_{U, K}(x) \leq c_{d, \phi}^{t}\,t^{5\,t + c_{1, \phi}}\, P^{d}\, (\log P)^{2\,s + c_{2, \phi}}\, D_{K}\,(\log D_{K})^{2d-2}, \label{borne.g=0}
\end{equation}
o\`u $c_{1, \phi}$ et $c_{2, \phi}$ ne d\'ependent que de $\phi$ et $c_{d, \phi}$ d\'epend en plus, de $d$.
\end{remarque}

\`a ce sujet, et contrairement \`a ce qui a \'et\'e d\'emontr\'e ici, D. Poulakis obtient des r\'esultats explicites. 
Soit $F \in K[X,Y]$  un polyn\^ome absolument irr\'eductible d\'efinissant une courbe de genre nul ayant au moins trois points \`a l'infini. Notons $N$ son degr\'e et $H_{K}(F)$ son hauteur multiplicative. Soient  $x$ et $y$ des $T$-entiers v\'erifiant $F(x,y)=0$. Dans \cite{poulakis.colloq}, il montre que
$$\max\{h_{K}(x),h_{K}(y)\} \leq   \gamma_{d, N, t}\, H_{K}(F)^{6000\,N^{3N+10} \,d^{2}\,t} \, \mathcal{P}^{300\,N^{3N+4} \,d\,t }\, D_{K}^{65\,N^{3N+4}\,d\,t},      $$
o\`u $\mathcal{P}= \max_{\vp \in S}N(\vp)$ et $\gamma = N^{10^{6}\,N^{5N+10}\,d^{2}\,t^{3}}$. On a  $P\leq \mathcal{P}\leq P^{d}$. 

La borne (\ref{borne.g=0}) donne une meilleure d\'ependance en le discriminant $D_{K}$; en particulier, l'exposant de $D_{K}$ ne d\'epend pas en l'ensemble $S$. Elle donne aussi une meilleure d\'ependance en l'ensemble de places $S$ quand son cardinal est petit par rapport \`a $P$.

Dans le cas particulier o\`u $T=M_{K}^{\infty}$,  D. Poulakis \cite{poulakis.hung} obtient :
$$\max\{h_{K}(x),h_{K}(y)\} \leq \xi_{d,N}\, D_{K}^{4730\,N^{9}}\, H_{K}(F)^{10^{9}\,N^{35}},$$
o\`u $\xi_{d,N} = d^{17 d\, N^{3}}(9N^{5N+4})^{10^{9}\,d\,N^{35}}$. On remarque ici en particulier que l'exposant de $D_{K}$ ne d\'epend pas du degr\'e $d$.

\medskip


\bigskip

\subsection{Quelques \'equations diophantiennes.}\label{eq.dio.classiques}


\begin{definition}\label{type.baker}
On dira qu'une courbe alg\'ebrique $U$ est contr\^ol\'ee par une autre courbe $V$ du point de vue des points entiers, si la finitude des points ($S$-)entiers de $U$ peut \^etre d\'eduite de celle des points ($S$-)entiers de $V$ \`a l'aide de morphismes comme ceux des propositions \ref{morph.fini} ou \ref{revet.etale}.

\end{definition}

En particulier,  $U$ est contr\^ol\'ee par $\droitemoins$ s'il existe $m\geq 1$ et un diagramme


\begin{displaymath}
\xymatrix@C=1cm@R=.5cm{
   & U_{1} \ar[dl]^{\psi_{1}} \ar[dr]^{\phi_{1}} &     &  U_{2} \ar[dl]^{\psi_{2}} \ar@{.}[dr]  &   & U_{m}  \ar[dl]^{\psi_{m}} \ar[dr]^{\phi_{m}}&   \\
U=V_{0} &                                        &  V_{1} &                                  &  V_{m-1} &   &  V_{m} = \droitemoins \\
           }
\end{displaymath}
avec, pour tout $i \in \{1, \ldots , m\}$, le morphisme $\phi_{i}$ fini et le morphisme $\psi_{i}$ fini, surjectif et non ramifi\'e en dehors des \g points \`a l'infini"  de $V_{i-1}$.

Voici quelques exemples de courbes contr\^ol\'ees par $\droitemoins$ : les courbes donn\'ees par une \'equation de Thue, les elliptiques, les superelliptiques et les courbes hyperelliptiques ayant un point de Weierstrass \`a l'infini.


\begin{cor}\label{thm.baker-serre}

Soit $U$ une courbe affine d\'efinie sur un corps de nombres $K$ de degr\'e $d =[K:\rat]$ telle que $\chi(U) < 0$ et $h_U$ une hauteur  sur $U$. Supposons que $U$ est contr\^ol\'ee par $\droitemoins$ au sens de la d\'efinition \ref{type.baker}. Soit $S$ un ensemble fini de places de $K$ et $P$ le maximum de ses caract\'eristiques r\'esiduelles.
Pour tout point $S$-entier $x$ de $U$, on a
$$h_U(x) \leq  k_{d}^{\card (S)}\, \card (S)^{k_{1}\,\card (S) + k_{2}}\, P^{k_{3}\,d} \,(\log P)^{k_{4}\,\card (S) + k_{5}} \,e^{\gamma_{d}} \,\gamma_{d}^{k_{6}\,d -2},$$   
o\`u $\gamma_{d} = k_{7}\,(\log D_{K} + \Sigma_{S} + k_{8} + d\,k_{9} \, \frac{\Sigma_{S} + k_{10}}{\log(\Sigma_{S}) + k_{11}})$,  les r\'eels $k_{j}$ d\'ependent de $U$ et $k_{d}$ d\'epend en plus de $d$.
\end{cor}

Y. Bugeaud \cite{bugeaud.super} donne une borne pour la hauteur des points $S$-entiers des courbes superelliptiques. Sa borne a l'avantage, par rapport \`a celle du corollaire \ref{thm.baker-serre},  de rendre explicites les exposants qu'y apparaîssent. L'ordre de grandeur des deux bornes est le m\^eme en ce qui concerne la d\'ependance en le discriminant du corps $K$, et aussi en l'ensemble de places $S$, si l'on fixe son cardinal.  On pouvait pr\'evoir cette ressemblance des bornes puisque les deux r\'esultats sont d\'eduits de \cite{bugeaud-gyory}.

\medskip

\noindent{\bf D\'emonstration du corollaire \ref{thm.baker-serre}.}

Nous appliquons  les propositions \ref{morph.fini} et \ref{revet.etale}, respectivement  aux morphismes $\phi_{i}$ et $\psi_{i}$ qui lient $U$ \`a $\droitemoins$ (en suivant le chapitre 8.4 de \cite{serre}, par exemple), et \`a la fonction $B_{\droitemoins}(u, v, w, z, d) = c_{d}^{v}\,v^{5(\,v + 2)}\,w^{d}  \,(\log w )^{2\,v} \, e^{z}\,z^{2\,d -2}$ donn\'ee par le corollaire \ref{cor.P1}. Si le morphisme $\phi : X \longrightarrow \droitemoins$ est fini, d'apr\`es la proposition \ref{morph.fini}, la courbe $X$ v\'erifie la propri\'et\'e \ref{prop.finitude} pour la fonction $B_{X}(u, v, w, z , d) = c_{3,d}^{v}\,v^{5\,v + c_{1}}\, w^{d}\, (\log w)^{2\,v + c_{2}}\, e^{z}\,z^{2d-2}$, o\`u $c_{1}$ et $c_{2}$ d\'ependent du morphisme $\phi$  et $c_{3,d}$ d\'epend en plus de $d$. Si $\psi : X \longrightarrow Y$ est un rev\^etement \'etale, d'apr\'es la proposition \ref{revet.etale}, la courbe $Y$ v\'erifie la propri\'et\'e \ref{prop.finitude} pour la fonction $B_{Y}(u, v, w, z, d) = c_{4, d}^{v}\, v^{c_{5}\,v + c_{6}}\, w^{c_{7}\,d} \,(\log w)^{c_{8}\,v + c_{9}} \,e^{\gamma_{d}} \,\gamma_{d}^{c_{10}\,d -2}$, o\`u $\gamma_{d} = c_{11}\,(z + u + c_{12} + c_{0}\, d\,\log(c_{11}) \, \frac{u + c_{12}}{\log(u + c_{12})})$, et les r\'eels $c_{j}$ d\'ependent de $\psi$ et $c_{4,d}$ d\'epend en plus de $d$. On remarque que par des appliquations succesives des propositions \ref{morph.fini} et \ref{revet.etale} on ne change pas l'odre de grandeur de la borne de la hauteur des points entiers, m\^eme si les constantes $k_{j}$ sont modifi\'ees. 
\hfill $\Box$

\bigskip

En particulier, dans le cas d'une courbe elliptique nous obtenons le r\'esultat ci-dessous.

\begin{cor}\label{thm.ellip}
Soient $E$ une courbe affine d\'efinie sur un corps de nombres $K$ de degr\'e $d$ dont la compl\'et\'ee $\Ebarre$ est une courbe elliptique, et $h_{E}$ une hauteur d\'efinie sur la courbe $E$. Soit $S$ un ensemble fini de places de $K$ et $P$ le maximum de ses caract\'eristiques r\'esiduelles.
Pour tout point $S$-entier  $p$ de $E$, on a 
$$h_{E}(p)\leq  \gamma_{E}\,c_{d}^{s + c_{1, E}}\,s^{20\,s + c_{2, E}}\,P^{4\,d}\,(\log P)^{8\,s + c_{3,E}}\, e^{\gamma_{d}}\,\gamma_{d}^{8d-2},$$
o\`u $\gamma_{d} = 4\,(\log D_{K} + \Sigma_{S} + c_{4,E} + c_{0}\, \log 4 \, d\,\frac{\Sigma_{S} + c_{5,E}}{\log(\Sigma_{S} + c_{5, E})})$, les nombres $c_{i, E}$ d\'ependent de $E$ et $\gamma_{E}$ d\'epend du choix des morphismes liant $E$ \`a $\droitemoins$, ainsi que de celui de la hauteur.

\end{cor}

\noindent{\bf D\'emonstration du corollaire \ref{thm.ellip}.}
Supposons que la courbe $E$ a un unique point \g \`a l'infini" $O$ qui est $K$-rationnel et prenons-le comme origine de sa loi de groupe. La multiplication par $2$, not\'ee $\psi$, est \'etale, de degr\'e $4$, donc l'image inverse de l'origine $O$ consiste en quatre points : $O, E_{1}, E_{2}, E_{3}$. 
Si notre courbe $E$ est donn\'ee par l'\'equation de Weierstrass $y^{2}=f(x)$ o\`u $f(x)= (x-e_{1})(x-e_{2})(x-e_{3})$, alors $E_{i}=(e_{i}:1)$.
En composant le morphisme donn\'e par la $x$-coordonn\'ee avec le morphisme qui envoie $(r:1)$ sur $((r-e_{1})(e_{2}-e_{3}) : (e_{2}-e_{1})(r-e_{3}))$, nous obtenons un morphisme fini $\phi : E' = \Ebarre \setminus \{O, E_{1}, E_{2}, E_{3}\} \vers \espproj^{1}\setminus\{0,1,\infty\}$, de degr\'e 2. 

D'apr\`es le corollaire \ref{cor.P1}, $\droitemoins$ v\'erifie la propri\'et\'e \ref{prop.finitude} pour la fonction 
$$B_{\droitemoins}(u, v, w, z, d) = c_{d}^{v}\,v^{5(\,v + 2)}\,w^{d}  \,(\log w )^{2\,v} \, e^{z}\,z^{2\,d -2}.$$

 En appliquant la proposition \ref{morph.fini} au morphisme $\phi$, nous en d\'eduisons que la courbe $E'$ v\'erifie la propri\'et\'e \ref{prop.finitude} pour la fonction 
$$B_{E'}(u,v,w,z,d) = c_{\phi}\,c_{d}^{v + v_{\phi}}\,(v+v_{\phi})^{5(v+v_{\phi}+2)}\,(w + w_{\phi})^{d}\,(\log(w+w_{\phi}))^{2\,(v + v_{\phi})}\, e^{z}\,z^{2\,d-2},$$
o\`u $c_{\phi}$ ne d\'epend que de $\phi$, $c_{d}$ que de $d$, $v_{\phi} = \card (S_{\phi})$, $w_{\phi} = \max P(S_{\phi})$ et $S_{\phi}$ est l'ensemble des id\'eaux premiers en lesquels $\phi$ ou $E'$ ont \g mauvaise r\'eduction".

En appliquant la proposition \ref{revet.etale} au morphisme $\psi$, nous en d\'eduisons que $E$ v\'erifie la propri\'et\'e de finitude \ref{prop.finitude} pour la fonction  
$$B_{E}(u,v,w,z,d)=c_{\psi}\,c_{\phi}\,c_{d}^{4\,(v+v_{\psi})+v_{\phi}}\,(4\,(v+v_{\psi})+ v_{\phi})^{5\,(4\,(v+v_{\psi})+ v_{\phi}) +2}\,(w+w_{\psi}+w_{\phi})^{4\,d}\times $$
$$\times (\log(w+w_{\psi}+w_{\phi}))^{2\,(4\,(v+v_{\psi})+v_{\phi})}\,e^{\gamma_{d}}\,\gamma_{d}^{8\,d-2},$$
o\`u $\gamma_{d} = 4\,(z + u + u_{\psi} + c_{0}\,\log 4 \,d\,\frac{u+u_{\psi}}{\log(u+u_{\psi})})$, $c_{\psi}$ ne d\'epend que de $\psi$, $c_{d}$ que de $d$, $u_{\psi} = \Sigma_{S_{\psi}}$, $v_{\psi} = \card (S_{\psi})$, $w_{\psi} = \max P(S_{\psi})$ et $S_{\psi}$ est l'ensemble des id\'eaux premiers en lesquels $\psi, E'$ ou $E$ ont \g mauvaise r\'eduction". Ceci d\'emontre le r\'esultat.

\hfill $\Box$

\medskip

 Posons $f(x)= x^{3} + a x + b$, avec $a$ et $b$ dans $O_{K,S}$ et tels que $4a^{3} + 27 b^{2} \ne 0$. Soient $(x,y)$ dans $O_{K,S}^{2}$ et v\'erifiant $y^{2}= f(x)$. Pour   $K = \rat$,  L. Hadju et T. Herendi \cite{hadju-herendi} montrent  que :
$$\max \{h(x), h(y)\} \leq (c_{1,f} \,s + c_{2,f})\, 10^{38\, s + 86}\, (s + 1)^{20\,s + 35} \,P^{24}\, (\log^{+}(P))^{4\,s +2},$$
o\`u $c_{1,f}$ et $c_{2,f}$ ne d\'ependent que en le polyn\^ome $f$.

 Quand l'ensemble de places $S$ est r\'eduit \`a l'ensemble des places archim\'ediennes, on a le r\'esultat de Y. Bugeaud \cite{bugeaud.ellip} :
 $$\max \{h(x),h(y)\} \leq c_{d, f}\, D_{K}^{6}\,(\log D_{K})^{12\,d +1},$$
o\`u $c_{d, f}$ d\'epend uniquement en le degr\'e $d$ de $K$ et en le polyn\^ome $f$; alors que le corollaire \ref{thm.ellip} donne une borne l\'eg\`erement meilleure : $c_{E,d, S} \,D_{K}^{4}\,(\log D_{K})^{8\,d-2}$.

\bigskip

%

\textbf{Remerciements.} L'essentiel de ce travail est inclus dans ma th\`ese de doctorat, r\'ealis\'ee \`a l'Institut de Math\'ematiques de Jussieu. Je tiens \`a remercier mes directeurs de th\`ese, Marc Hindry et Michel Waldschmidt, ainsi que Joseph Oesterl\'e, pour de nombreuses discussions sur ce sujet.

%
%

\bigskip
\bigskip

\begin{quote}
Andrea Surroca Ortiz\\
\'Ecole Polytechnique F\'ed\'erale de Lausanne\\
FSB IMB CSAG\\
MA C3 635 (B\^at. MA)\\
Station 8\\
CH-1015 Lausanne\\
andrea.surroca@epfl.ch
\end{quote}


\begin{thebibliography}{Dillo99}



\bibitem[1]{belyi} \Belyi, G.V., {\em On Galois extensions of a maximal cyclotomic field}, Math. USSR Izveztija, Vol. 14, n°2 (1980).

\bibitem[2]{bilu} Bilu, Yu. F., {\em Quantitative Siegel's theorem for Galois coverings}, Compositio Math. 106 (1997), no. 2, 125-158.

\bibitem[3]{bugeaud.ellip} Bugeaud, Y. {\em On the size of integer solutions of elliptic equations}, Bull. Austral. Math.Soc. vol. 57 (1998), 199-206.

\bibitem[4]{bugeaud.super} Bugeaud, Y. {\em Bounds for the solutions of superelliptic equations}, Compositio Math. 107 (1997), 187-219.

\bibitem[5]{bugeaud-gyory} Bugeaud, Y., Gy\"ory, K., {\em Bounds for the solutions of unit equations}, Acta Arith. 74 (1996), 273-292.

\bibitem[6]{elkies} Elkies, N.D., {\em ABC implies Mordell}, Int. Math. Res. Not. 7 (1991), 99-109. 


\bibitem[7]{esnault-viehweg} Esnault, H., Viehweg, E., {\em Effective bounds for semipositive sheaves and for the height of points of curves over complex function fields}, Compos. Math. 76, No.1/2 (1990), 69-85.


\bibitem[8]{franken} Frankenhuysen, M. van, {\em The ABC conjecture implies Roth's theorem and Mordell's conjecture}, Matem\'atica Contempor\^anea, Vol. 16 (1999), 45-72. 


\bibitem[9]{sga} Grothendieck, A., {\em Rev\^etements \'etales et groupe fondamental}, S\'eminaire de g\'eom\'etrie alg\'ebrique du Bois Marie 1960/61 (SGA 1), Lecture Notes in Mathematics, 224, Springer-Verlag, Berlin-Heidelberg-New York, (1971).


\bibitem[10]{hadju-herendi}  Hadju, L., Herendi, T.,  {\em Explicit bounds for the Solutions of Elliptic Equations with Rational Coefficients}, J. symbolic Computation 25 (1998), 361-366.

\bibitem[11]{hindry-silverman.inv}  Hindry, M.,  Silverman, J.H., {\em The canonical height and integral points on elliptic curves}, Invent. Math. 93, No.2 (1988), 419-450.


\bibitem[12]{hindry-silverman}  Hindry, M.,  Silverman, J.H., {\em Diophantine
    Geometry. An introduction}, Graduate texts in mathematics,
    Springer-Verlag, New York, 2000.

\bibitem[13]{kubert-lang} Kubert, D., Lang, S., {\em Units in the Modular Function Field  I}, Math. Ann. 218 (1975), 67-96.


\bibitem[14]{mason.dio.equa} Mason, R.C., {\em Diophantine Equations over Function Fields}, London Mathematical Society Lecture Note Series 96, Cambridge University Press, Cambridge, 1984.


\bibitem[15]{masser} Masser, D.W., {\em On $abc$ and discriminants}, Proceedings of the American Mathematical Society, Volume 130, Number 11, (2002), 3141-3150.
 

\bibitem[16]{moret-bailly} Moret-Bailly, L., {\em Hauteurs et classes de Chern sur les surfaces arithm\'etiques}, Ast\'erisque 183, (1990), 37-58.

\bibitem[17]{oesterle} Oesterl\'e, J., {\em Nouvelles approches du ``th\'eor\`eme'' de Fermat}, S\'eminaire Bourbaki, 40\`eme ann\'eee, n°694, f\'evrier 1998.


\bibitem[18]{poulakis.hung} Poulakis, D., {\em Bounds for the size of integral points on curves of genus zero}, Acta Math. Hung. 93, n°4, 327-346 (2001).

\bibitem[19]{poulakis.colloq} Poulakis, D., {\em Points entiers sur les courbes de genre 0}, Colloq. Math. 66 (1993), n°1, 1-7.

\bibitem[20]{serre} Serre, J.-P., {\em Lectures on the Mordell-Weil theorem}. Third edition. Aspects of Mathematics. Vieweg (1997).

\bibitem[21]{serre-chebotarev} Serre, J.-P., {\em Quelques appliquations du th\'eor\`eme de densit\'e de Chebotarev}, Inst. Hautes \'Etudes Sci. Publ. Math. n°54, (1981), 323-401. 


\bibitem[22]{siegel} Siegel, C.L., {\em \"Uber einige Anwendungen diophantischer Approximationen}, Abh. Pr. Akad. Wiss. 1 (1929) 41-69 (Ges. Abh., I, 209-266).  


\bibitem[23]{stewart-tijdeman} Stewart, C.L., Tijdeman R., {\em On the Oesterl\'e-Masser Conjecture}, Monatsh. Math., 102, (1986), 251-257. 

\bibitem[24]{stewart-yu} Stewart, C.L., Yu, Kunrui, {\em On the $abc$ conjecture, II}, Duke Math. J. 108 (2001), n°1, 169-181.  

\bibitem[25]{mathese} Surroca, A., {\em M\'ethodes de transcendance et g\'eom\'etrie diophantienne}, Th\`ese de doctorat de l'Universit\'e Paris VI, soutenue le 1er d\'ecembre 2003. (100 pages) http://www.alg-geo.epfl.ch/$\sim$surroca

\bibitem[26]{parma} Surroca, A., {\em Siegel's theorem and the $abc$ conjecture}, Proceedings of the Secondo Convegno Italiano di Teoria dei Numeri, Parma, Nov. 2003. Riv. Mat. Univ. Parma (7) (2004), 323-332.


\bibitem[27]{szpiro} Szpiro, L., {\em Propri\'et\'es num\'eriques du faisceau dualisant r\'elatif}, Seminaire sur les pinceaux de courbes de genre au moins deux, Ast\'erisque 86 (1981), 44-78.


\bibitem[28]{vojta} Vojta, P., {\em On algebraic points on curves}, Compos. math. 78 (1991), no. 1, 29-36.



\end{thebibliography}
\end{document}